\documentclass[11pt]{article}

\usepackage{amsmath}
\usepackage{amssymb,amsbsy}
\usepackage{graphicx}
\usepackage{dsfont}
\usepackage{enumerate}
\usepackage[margin=1in]{geometry}
\usepackage{tikz}
\usetikzlibrary{arrows,decorations.pathmorphing,backgrounds,positioning,fit,petri,mindmap,shapes.geometric,decorations.pathreplacing}

\usepackage{pdfsync}
\usepackage{hyperref}

\allowdisplaybreaks[1]

\newtheorem{theorem}{Theorem}

\newtheorem{lemma}{Lemma}

\newtheorem{conjecture}{Conjecture}

\renewcommand{\phi}{\varphi}

\renewcommand{\P}{\mathbb{P}}
\newcommand{\E}{\mathbb{E}}

\newcommand{\R}{\mathbb{R}}

\def\ds1{\mathds{1}}
\renewcommand{\epsilon}{\varepsilon}
\newcommand{\eps}{\epsilon}
\newcommand\Var{{\dsV\text{ar}}\,}
\newcommand\dsV{\mathbb{V}}

\newcommand{\wh}{\widehat}

\renewcommand{\tilde}{\widetilde}

\newlength{\minipagewidth}
\newcommand{\beq}{\begin{equation}}
\newcommand{\eeq}{\end{equation}}

\newcommand{\beqa}{\begin{eqnarray}}
\newcommand{\eeqa}{\end{eqnarray}}

\newcommand{\beqan}{\begin{eqnarray*}}
\newcommand{\eeqan}{\end{eqnarray*}}

\def\ba#1\ea{\begin{align*}#1\end{align*}} 
\def\banum#1\eanum{\begin{align}#1\end{align}} 

\newcommand{\TV}{\mathrm{TV}}
\newcommand{\Tr}{\mathrm{Tr}} 

\newcommand{\EE}{\mathbb{E}}
\newcommand{\RR}{\mathbb{R}}

\newcommand{\Proj}{\mathrm{Proj}}
\newcommand{\BlackBox}{\rule{1.5ex}{1.5ex}}  
\newenvironment{proof}{\par\noindent{\bf Proof\ }}{\hfill\BlackBox\\[2mm]}

\begin{document}
\title{Testing for high-dimensional geometry in random graphs}
\author{
	S{\'e}bastien Bubeck
	\thanks{Microsoft Research and Princeton University; \texttt{sebubeck@microsoft.com}.}
	\and
	Jian Ding
	\thanks{University of Chicago; \texttt{jianding@galton.uchicago.edu}.}
	\and 
	Ronen Eldan 
	\thanks{University of Washington; \texttt{roneneldan@gmail.com}.}
	\and 
	Mikl\'os Z.\ R\'acz
	\thanks{University of California, Berkeley; \texttt{racz@stat.berkeley.edu}.}
}
\date{\today}

\maketitle

\begin{abstract}
We study the problem of detecting the presence of an underlying high-dimensional geometric structure in a random graph. Under the null hypothesis, the observed graph is a realization of an Erd{\H{o}}s-R{\'e}nyi random graph $G(n,p)$. Under the alternative, the graph is generated from the $G(n,p,d)$ model, where each vertex corresponds to a latent independent random vector uniformly distributed on the sphere $\mathbb{S}^{d-1}$, and two vertices are connected if the corresponding latent vectors are close enough. In the dense regime (i.e., $p$ is a constant), we propose a near-optimal and computationally efficient testing procedure based on a new quantity which we call signed triangles. The proof of the detection lower bound is based on a new bound on the total variation distance between a Wishart matrix and an appropriately normalized GOE matrix. In the sparse regime, we make a conjecture for the optimal detection boundary. We conclude the paper with some preliminary steps on the problem of estimating the 
dimension in $G(n,p,d)$.
\end{abstract}

\section{Introduction}
Extracting information from large graphs has become an important statistical problem, as network data is now common in various fields. Whether one talks about social networks, gene networks, or (biological) neural networks, in all cases there is a rapidly growing industry dedicated to extracting knowledge from these graphs. A particularly important task in this spirit is to learn a useful {\em representation} of the vertices, that is, a mapping from the vertices to some metric space (usually $\R^d$). 
In this paper, we take a step back and study the hypothesis testing problem that underlies these investigations: 
given a large graph, one wants to tell if the observed connections result from a latent geometrical structure in the vertices, or if they are purely random.

\medskip

The null hypothesis that we consider is that the observed graph $G$ on $n$ vertices has been generated by the standard Erd\H{o}s-R\'enyi random graph $G(n,p)$ \cite{ER60}, where each edge appears independently with probability $p$:
\[
H_0 : G \sim G(n,p) .
\] 
On the other hand, for the alternative, we consider the simplest model of a random geometric graph---recall that a geometric graph is such that each vertex is labeled with a point in some metric space, and an edge is present between two vertices if the distance between the corresponding labels is smaller than some prespecified threshold. We focus on the case where the underlying metric space is the Euclidean sphere $\mathbb{S}^{d-1}=\{x\in \R^d: \|x\|_2=1\}$, and the latent labels are i.i.d. uniform random vectors in $\mathbb{S}^{d-1}$. We denote this model by $G(n,p,d)$, where $p$ is the probability of an edge between two vertices ($p$ determines the threshold distance for connection), and thus
$$H_1 : G \sim G(n, p, d) .$$
Slightly more formally, $G(n,p,d)$ is defined as follows. 
Let $X_1,\ldots,X_n$ be independent
random vectors, uniformly distributed on $\mathbb{S}^{d-1}$. In $G(n,p,d)$, distinct vertices $i \in [n]$ and $j \in [n]$ are connected by an edge if and only if $\langle X_i,X_j \rangle \ge t_{p,d}$, where the threshold
value $t_{p,d} \in [-1,1]$ is such that $\P\left( \langle X_1,X_2 \rangle \ge t_{p,d} \right) = p$. 

\medskip

In contrast to most previous works (see Section~\ref{sec:relatedwork}), the focus of this paper is on the {\em high-dimensional} situation, where $d$ can scale as a function of $n$. This point of view is in line with recent advances in all areas of applied mathematics, and in particular statistics and learning theory, where high-dimensional feature spaces are becoming the new norm. We also consider both the dense regime, where $p$ is a constant independent of $n$, and the sparse regime, where $p=c/n$ for some constant $c$.

\medskip

A natural test to uncover geometric structure is to count the number of triangles in $G$. Indeed, in a purely random scenario, vertex $u$ being connected to both $v$ and $w$ says nothing about whether $v$ and $w$ are connected. On the other hand, in a geometric setting this would imply that $v$ and $w$ are close to each other due to the triangle inequality, thus increasing the probability of a connection between them. 
This, in turn, implies that the expected number of triangles is larger in the geometric setting, given the same edge density. More broadly, the number of triangles is a basic statistic in the analysis of networks, as it gives an indication of the dependencies in the networks, see, e.g.,~\cite{GRS14} and the references therein. 

\medskip

One of the main results of this paper is that in the dense regime, there are much more powerful statistics than triangles to distinguish between $G(n,p)$ and $G(n,p,d)$. In particular, we propose a near-optimal statistic based on a new quantity which we call {\em signed triangles}. 
To put 
it succintly, 
if $A$ denotes the adjacency matrix of a graph $G$, 
then the total number of triangles is equal to $\mathrm{Tr}(A^3)$, 
while the total ``number" of signed triangles is defined as $\mathrm{Tr} \left( \left( A - p \left( J - I \right) \right)^3 \right)$, 
where $I$ is the identity matrix, and $J$ is the matrix with every entry equal to $1$. 
The key insight motivating this definition is that the variance of signed triangles is dramatically smaller than the variance of triangles, due to the cancellations introduced by the centering of the adjacency matrix. 
While being elementary, this idea of using cancellations of various patterns seems to be new in the network analysis literature. It could have applications beyond the problem studied in this paper, in particular because it is quite different from other popular statistics, such as those coming from the method of moments \cite{BCL11} (where patterns are always counted positvely), or those with the flavor of the clustering coefficient \cite{WS98} (where ratios of counts are considered). 
Indeed, the key aspect of signed triangles is that different patterns are reweighted with positive and negative weights: 
triangles and induced single edges are counted positively, while induced wedges and independent sets on three vertices are counted negatively.

\subsection{Related work} \label{sec:relatedwork}
In the social networks literature the latent metric space is referred to as the {\em social space}. Starting with \cite{HRH02}, there have been numerous works on estimating positions in the social space for individuals in a social network. Various models are considered in this literature, and interestingly some of them are very closely related to $G(n,p,d)$, see~\cite[Section~2.2]{HRH02}. Most papers focus on various approximations to the maximum likelihood estimator, and some theoretical results have been obtained with other methods such as counting the number of common neighbors, see \cite{SCM10, ACKS13}. 
In the present paper we settle for a less ambitious goal, as we are not learning the representation but simply testing for the presence of a meaningful representation. 
On the other hand, we  obtain much more precise results, such as an almost tight characterization of couples $(n,d)$ for which testing between $G(n,p)$ and $G(n,p,d)$ is possible in the dense regime. Most importantly, this provides, to the best of our knowledge, the first result for the high-dimensional setting, as all previous works on social space inference were in the low-dimensional regime (i.e., $d$ is fixed and $n \to \infty$).

\medskip

In probability theory, models such as $G(n,p,d)$ have been studied for a long time in the low-dimensional regime, see, e.g.,~\cite{Pen03}.  The high-dimensional setting was first investigated recently in~\cite{DGLU11}. In this paper it was observed that with $n$ fixed and $d \to \infty$, $G(n,p,d)$ converges in total variation to $G(n,p)$. 
In other words, if the dimension $d$ is very large compared to $n$, then one cannot distinguish between $G(n,p)$ and $G(n,p,d)$. 
Our main result is that, in fact, in the dense regime, there is a phase transition when $d$ is of order $n^3$, with the total variation distance going from $1$ to $0$. 
In spite of previous works, the high-dimensional $G(n,p,d)$ remains mysterious in many ways. Essentially, in the dense regime, the only graph parameter which is well understood is the clique number, due to the results of~\cite{DGLU11, ABL13}, while in the sparse case basically nothing is known. 
One of the technical contributions of the present paper is to compute rather precisely the probability of a triangle in the sparse case.

\medskip

This paper can be seen as part of a series of recent papers studying hypothesis testing in random graph models~\cite{AV13, AV13b, BMR14}. For instance, Arias-Castro and Verzelen consider the problem of testing for {\em community structure}, while we test for {\em geometric structure}. More precisely, their null is identical to ours, that is $G \sim G(n,p)$, while for the alternative they consider the model $G(n,p,k,q)$ which differs from $G(n,p)$ by having a random subset of $k$ vertices for which edges appear between them with probability $q$. This problem is closely related to community detection in stochastic block models, a problem which has recently attracted a lot of attention, see, e.g.,~\cite{NMS14,ABH14} and the references therein. Interestingly, when testing for community structure, one of the main obstacles is computational rather than information-theoretic. For example, it is obvious that one can distinguish between $G(n,1/2)$ and $G(n,1/2,k,1)$ as long as $k\gg \log(n)$, though when $k=o(\sqrt{n}
)$, no 
polynomial test is known for 
this problem (referred to as the {\em planted clique problem}), see, e.g.,~\cite{AKS98, BR13}. On the contrary, in the context of testing for geometrical structure, we show that polynomial time methods are near-optimal, at least for the dense regime, since one can efficiently compute signed triangles. This is perhaps surprising, as the worst-case version of our problem, namely recognizing if a graph can be realized as a geometric graph, is known to be NP-hard~\cite{BK98}. 
Finally, we also note that our new signed triangles statistic is closely related to the tensor introduced in~\cite{FK08} for the planted clique problem. While~\cite{FK08} computes the spectral norm of this tensor, here we simply sum its entries (in the case $p=1/2$).

\subsection{Contributions and content of the paper}
The main objective of the paper is to identify the boundary of testability between $G(n,p)$ and $G(n,p,d)$, in both the dense and the sparse regimes. To put it differently, we are interested in studying the total variation distance between these two models, denoted by $\TV\left(G(n,p), G(n,p,d)\right)$. Recall that in~\cite{DGLU11} it was proved that, if $n$ is fixed and $d \to \infty$, then
\[
\TV\left(G(n,p), G(n,p,d)\right) \to 0.\footnote{More precisely, they show that this convergence happens when $d$ is exponential in $n$.}
\]
Given this result, our question of interest becomes the following: how large can $d$ be, as a function of $n$, so that $\TV\left(G(n,p), G(n,p,d)\right)$ remains bounded away from $0$ (or even becomes close to $1$)? The core of our contribution to this question is summarized by the following theorem.

\begin{theorem} \label{th:mainresult}
\begin{enumerate}[(a)]
 \item\label{th:main_denseUB} Let $p \in (0,1)$ be fixed, and assume that $d/n^3 \to 0$. Then 
\begin{equation} \label{eq:mainresult1}
\TV\left(G(n,p), G(n,p,d) \right) \to 1.
\end{equation}
 \item\label{th:main_sparseUB} Let $c > 0$ be fixed, and assume that $d/\log^3(n) \to 0$. Then 
\begin{equation} \label{eq:mainresult2}
\TV\left(G\left(n,\frac{c}{n}\right), G\left(n,\frac{c}{n},d\right) \right) \to 1.
\end{equation}
 \item\label{th:main_denseLB} Furthermore, if $d/n^3 \to \infty$, then 
\begin{equation} \label{eq:mainresult3}
\sup_{p \in \left[0,1\right]} \TV\left(G(n,p), G(n,p,d)\right) \to 0 .
\end{equation}
\end{enumerate}
\end{theorem}
The results of Theorem~\ref{th:mainresult}\eqref{th:main_denseUB} and~\ref{th:mainresult}\eqref{th:main_denseLB} 
tightly characterize the dense regime. 
We conjecture that the result of Theorem~\ref{th:mainresult}\eqref{th:main_sparseUB} 
is tight for the sparse regime:
\begin{conjecture} \label{conj}
Let $c > 0$ be fixed, and assume that $d/\log^3(n) \to \infty$. Then 
$$\TV\left(G\left(n,\frac{c}{n}\right), G\left(n,\frac{c}{n},d\right) \right) \to 0.$$
\end{conjecture}
In the following we outline the methods we use to prove Theorem \ref{th:mainresult}.

\medskip

We first discuss the proof of part~\eqref{th:main_denseUB} of Theorem~\ref{th:mainresult}, 
where our main methodological contribution lies. As pointed out above, a natural test to consider is the number of triangles:
\begin{equation}\label{eq:triangle_def}
 T(G) := \sum_{\{i,j,k\} \in \binom{[n]}{3}} A_{i,j} A_{i,k} A_{j,k},
\end{equation}
where $A$ denotes the adjacency matrix of the graph $G$ with vertex set $\left[n\right]$, and so $A_{i,j} = 1$ if vertices $i$ and $j$ are connected by an edge, and $A_{i,j} = 0$ otherwise. 
An elementary calculation shows that $\E \ T(G(n,p)) = \binom{n}{3} p^3$, while $\Var(T(G(n,p)))$ is of order $n^4$. A key calculation of the paper is to show that $\E \ T(G(n,p,d)) \geq \binom{n}{3} p^3 (1 + C_p / \sqrt{d})$ for some constant $C_p$ that depends only on $p$, see Lemma~\ref{lem:probatriangle-p}. Intuitively, this shows that triangles have some power as long as 
\[
 \left|\E \ T(G(n,p)) - \E \ T(G(n,p,d)) \right| \gg \sqrt{\Var(T(G(n,p)))}, 
\]
which is equivalent to $d \ll n^2$. 
One of the main contributions of the paper is to introduce a statistic that is asymptotically powerful as long as $d \ll n^3$, which we refer to as the number of {\em signed triangles}:
\begin{equation}\label{eq:signed_triangle_def}
\tau(G) := \sum_{\{i,j,k\} \in \binom{[n]}{3}} (A_{i,j}-p) (A_{i,k}-p) (A_{j,k}-p) .
\end{equation}
The key point of signed triangles is the reduction of variance. 
Namely, $\Var(\tau(G(n,p)))$ is of order $n^3$, instead of $n^4$ for triangles. 
The improvement comes from the fact that 
$4$-vertex subgraphs with at least $5$ edges 
do not contribute to $\Var(\tau(G(n,p)))$, but they do contribute to $\Var(T(G(n,p)))$. 
We study $\tau$ in detail in Section~\ref{secdenseUB}, where we prove the following result (which implies Theorem~\ref{th:mainresult}\eqref{th:main_denseUB}).

\begin{theorem} \label{th:denseUB}
Let $p \in (0,1)$ be fixed, and assume that $d/n^3 \to 0$. Then 
\[
\TV\left(\tau(G(n,p)), \tau(G(n,p,d)) \right) \to 1.
\]
\end{theorem}

In the sparse regime we analyze the triangle statistic and prove the following theorem, which implies 
part~\eqref{th:main_sparseUB} of Theorem~\ref{th:mainresult}.

\begin{theorem} \label{sparseUB}
Let $c > 0$. If $d / \log^3(n) \to 0$, then 
\[
\TV\left(T\left(G\left(n,\frac{c}{n}\right)\right), T\left(G\left(n,\frac{c}{n},d\right)\right) \right) \to 1.
\]
\end{theorem}
In contrast with the dense regime, in the sparse regime the signed triangle statistic $\tau$ does not give significantly more power than the triangle statistic $T$. 
This is because in the sparse regime, with high probability, 
the graph does not contain any $4$-vertex subgraph with at least $5$ edges, 
which is where the improvement comes from in the dense regime. 

\medskip 

We believe that the upper bound $\log^3(n)$ in Theorem~\ref{sparseUB} cannot be significantly improved, as stated above in Conjecture~\ref{conj}. 
The main reason for this conjecture is that, when $d \gg \log^3(n)$, $G(n, c/n)$ and $G(n, c/n, d)$ seem to be locally equivalent (in particular, they both have the same Poisson number of triangles asymptotically). 
Thus the only way to distinguish between them would be to find an emergent global property which is significantly different under the two models, but this seems unlikely to exist.

\medskip

The final result of Theorem \ref{th:mainresult}, part~\eqref{th:main_denseLB}, complements part~\eqref{th:main_denseUB} by giving a matching lower bound. 
This bound is also valid for the sparse regime, though in this case we believe that it is not tight. 
The main idea behind the proof of this result is to view the random graphs in question as (essentially) the same function of appropriate random matrices. 
On the one hand, the random geometric graph $G\left(n,p,d\right)$ can be viewed as a function of an $n\times n$ Wishart matrix with $d$ degrees of freedom---i.e., a matrix of inner products of $n$ $d$-dimensional Gaussian vectors---denoted by $W(n,d)$. 
On the other hand, one can view $G\left( n, p \right)$ as a function of an $n \times n$ GOE random matrix---i.e., a symmetric matrix with i.i.d. Gaussian entries on and above the diagonal---denoted by $M(n)$. 
Moreover, these two functions are essentially the same. 
Theorem~\ref{th:mainresult}\eqref{th:main_denseLB} then follows from the following result on random matrices, stating that, if $d / n^3$ is very large, then the Wishart matrix has approximately the same law as an
appropriately centered and scaled GOE random matrix.
\begin{theorem}\label{thm:Wishart_GOE_intro}
Let $I_n$ denote the $n \times n$ identity matrix. If $d / n^3 \to \infty$, then
\begin{equation*}
  \TV \left( W \left( n, d \right), \sqrt{d} M \left( n \right) + d I_n \right) \to 0.
\end{equation*}
\end{theorem}
The random ensembles $W(n,d)$ and $M(n)$ are defined more precisely in Section~\ref{sec:denseLB}, where the above theorem is also proved.

\medskip

Finally, in Section \ref{sec:estimation} we touch upon the following question: for which values of $n$ and $d$ can the dimension $d$ be reconstructed by observing a sample of $G(n,p,d)$? We give the following bound for $p=1/2$, which can be considered as a proof of concept. 

\begin{theorem} \label{thm:estimation}
There exists a universal constant $C>0$, such that for all integers $n$ and $d_1<d_2$, one has
$$
\TV\left(G(n,1/2,d_1), G(n,1/2, d_2)\right) \geq 1 - C \left (\frac {d_1} n \right )^2.
$$
\end{theorem}

This bound is tight, as demonstrated by an older result of the third named author~\cite{Eldan11}, 
which states that when $d \gg n$, the graphs $G(n,1/2,d)$ and $G(n,1/2,d+1)$ are indistinguishable.

\begin{theorem}
There exists a universal constant $C>0$ such that for all integers $n<d$, 
\[
\TV\left(G(n,1/2,d), G(n,1/2, d+1)\right) \leq C \sqrt{ \left( \frac{d+1}{d-n} \right)^2 - 1}.
\]
\end{theorem}

\section{Estimates for the number of triangles in a geometric graph} \label{sec:triangle}

The point of this section is to give a lower bound for the expected number of triangles in the random geometric graph $G(n,p,d)$, using elementary estimates. 
To this end, let $X_1, X_2$, and $X_3$ be independent uniformly distributed points in $\mathbb{S}^{d-1}$. Consider the event 
\[
E = \Bigl \{\langle X_1, X_2 \rangle \geq t_{p, d}, \langle X_1, X_3 \rangle \geq t_{p, d}, \langle X_2, X_3 \rangle \geq t_{p, d} \Bigr \}
\]
that the corresponding vertices form a triangle; the expected number of triangles in $G(n,p,d)$ is thus $\E \ T(G(n,p,d)) = \binom{n}{3} \P(E)$. 
Our main result of this section is the following.

\begin{lemma} \label{lem:probatriangle-p}
There exists a universal constant $C > 0$ such that whenever $p<\tfrac 1 4$ we have that 
\begin{equation} \label{trianpsmall}
\P( E) \geq C^{-1} p^3 \frac{(\log \tfrac 1 p)^{3/2}}{\sqrt{d}}.
\end{equation}
Moreover, for every fixed $0<p<1$, there exists a constant $C_p > 0$ such that for all $d\geq C_p^{-1}$,
\begin{equation} \label{triang}
\P( E ) \geq p^3 \left (1 + \frac{C_p}{\sqrt{d}} \right ).
\end{equation}
\end{lemma}

Before we move on to the proof, we need some preliminary results, for which we need some additional notation. 
Denote by $\phi(x) = \left( 1/ \sqrt{2\pi} \right) e^{-x^2/2}$
the standard Gaussian density and by $\overline{\Phi}(x) = \int_{x}^\infty \phi(z) dz$ the associated complementary cumulative distribution function 
(note that $\overline{\Phi}$ is decreasing). 
Moreover, define
\begin{equation}\label{eq-density-Z}
f_{ d}(x) = \frac{\Gamma( d/2)}{\Gamma((d-1)/2) \sqrt{\pi}} (1 - x^2)^{(d-3)/2}, ~~ x \in \left[-1,1\right].
\end{equation}
A standard calculation gives that $f_d$ is the density of a one-dimensional marginal of a uniform random point on $\mathbb{S}^{d-1}$ (see~\cite[Section 2]{Sodin05}). 
Throughout the paper we use the notation $a \wedge b := \min \left\{ a, b \right\}$. 

\begin{lemma}\label{cor:tpd}
There exists a universal constant $C>0$ such that for all $0 < p \leq 1/2$ we have that
\begin{equation} \label{tpdest}
\left( C^{-1} (\tfrac 1 2 - p) \sqrt{\frac{\log (1/p)}{d}} \right) \wedge (1/2) \leq t_{p,d} \leq C \sqrt{\frac{\log(1/p)}{d}}.
\end{equation}
and
\begin{equation} \label{eq:fdtpd}
f_{d-1} (t_{p,d}) \geq C^{-1} d p t_{p, d}.
\end{equation}
\end{lemma}
\begin{proof}
We may assume that $d \geq 4$; the statement is easily checked when $d \leq 3$. 
We begin by recalling the following well-known inequalities (see, e.g.,~\cite[Chapter~6]{AbramowitzStegun:64}):
\begin{equation} \label{gammaineq}
\Gamma(x)\sqrt{x}/100 \leq \Gamma(x+1/2) \leq 2 \sqrt{x} \Gamma(x), ~~ \forall x > 1.
\end{equation} 
Therefore, we have for $0 < \theta < 1$ that
\begin{align}\label{eq-Psi-upper}
\int_{\theta}^{1} f_d(x) dx & \leq \frac{2 \sqrt{d}}{\sqrt{2\pi}} \int_\theta^1 (1-x^2)^{(d-3)/2} dx   = \frac{2 \sqrt{d}}{\sqrt{2\pi}} \int_0^{1-\theta} e^{\tfrac{d-3}{2} \log (1-(\theta + x)^2)} dx \nonumber \\
& \leq \frac{2 \sqrt{d}}{\sqrt{2\pi}} (1 - \theta^2)^{(d-3)/2} \int_0^{1-\theta} e^{\tfrac{d-3}{2} \left (\log (1-(\theta + x)^2) - \log (1-\theta^2) \right )} dx \nonumber \\
&\leq \frac{2 \sqrt{d}}{\sqrt{2\pi}} (1 - \theta^2)^{(d-3)/2} \int_0^1 e^{-\tfrac{d-3}{2} (\theta x + x^2) } d x  \nonumber \\
& \leq 8 \left (1  \wedge \frac{1}{\theta \sqrt{d}} \right )  (1 - \theta^2)^{(d-3)/2}\,.
\end{align}
Taking $\theta = t_{p,d}$ thus gives $p \leq 8 (1 - t_{p,d}^2)^{(d-3)/2}$ and so $\log (p/8) \leq - \tfrac{d-3}{2} t_{p,d}^2$, 
which implies the upper bound in \eqref{tpdest}. Moreover, the inequality in \eqref{eq-Psi-upper} also gives
$$ 
p \leq \frac{8}{t_{p,d} \sqrt{d}}  (1 - t_{p,d}^2)^{(d-3)/2} \stackrel{\eqref{gammaineq}}{\leq} \frac{8 C}{d t_{p,d}} f_d(t_{p,d})
$$
for a universal constant $C>0$, which proves \eqref{eq:fdtpd}.

\medskip

It remains to prove the lower bound in \eqref{tpdest}. First assume $10^{-6} \leq p \leq 1/2$.  For all $0\leq x\leq 1$ we see that $f_d(x) \leq 2 \sqrt{d}$, and so for all $0 \leq \theta \leq 1$ we have that 
$\int_0^\theta f_d(x ) dx \leq  2 \theta \sqrt{d} $.
Since $\int_0^{t_{p, d}} f_d(x )dx = 1/2 - p$, we get that $t_{p, d} \geq \frac{1/2-p}{2\sqrt{d}}$ as required. 

\medskip

Next, we consider the case for $0<p\leq 10^{-6}$. In this case it suffices to assume that $d$ is larger than an absolute constant so that the left hand side of \eqref{tpdest} is at most $1/2$.  By \eqref{gammaineq}, we see that for  $0 < \theta \leq 1/2$ we have that 
\begin{align}\label{eq-Psi-lower}
\int_{\theta}^{1} f_d(x) dx &\geq \frac{\sqrt{d}}{400} \int_\theta^1 (1-x^2)^{(d-3)/2} dx  \nonumber\\
&\geq \frac{ \sqrt{d}}{400}  \int_\theta^{\theta + \frac{1}{\sqrt{d}} \wedge \frac{1}{\theta d}} e^{-x^2 (d-3)} d x \geq 10^{-4} \left (1\wedge \frac{1}{\theta \sqrt{d}} \right ) (1 - \theta^2)^{(d-3)/2}\,.
\end{align}
Since $\int_{t_{p, d}}^1 f_d(x) dx = p$, the preceding inequality yields that $t_{p, d} \geq (c \sqrt{\log (1/p) /d}) \wedge (1/2)$ for some absolute constant $c>0$. This completes the proof.
\end{proof}

\medskip

We are finally ready to prove the main estimate of this section.

\medskip 

\begin{proof}\textbf{of Lemma~\ref{lem:probatriangle-p}} 
First, assume that $d \leq \tfrac 1 4 \log (1/p)$. 
Denote by $\angle(\cdot,\cdot)$ the geodesic distance on $\mathbb{S}^{d-1}$ and set $g(\theta) := \P(\angle(X_1, X_2) < \theta)$. 
A standard calculation shows that
$$
g(\theta) = \frac{(d-1) \pi^{(d-1)/2} }{ \Gamma\left (\tfrac{d+1}{2} \right ) }   \int_0^{\theta} \sin(x)^{d-2} dx.
$$
Using a change of variables and the fact that $\sin(x/2) \geq \sin(x) / 2$ for all $x \in \left[0, \pi \right]$, we have that 
\begin{equation} \label{eq:subdoubling}
g(\theta / 2) \geq 2^{-d} g(\theta), ~~ \forall \theta < \pi.
\end{equation} 
Next observe that one has, by the triangle inequality,
\begin{multline*}
E = \{\angle(X_1, X_2) < \arccos(t_{p,d}), \angle(X_1, X_3) < \arccos(t_{p,d}), \angle(X_2, X_3) < \arccos(t_{p,d})  \} \\
\supseteq \{ \angle(X_1, X_2) < \tfrac 1 2 \arccos(t_{p,d}),  \angle(X_1, X_3) < \tfrac 1 2 \arccos(t_{p,d}) \}.
\end{multline*}
Since $\angle(X_1,X_2)$ and $\angle(X_1,X_3)$ are independent, we get that
\[
\P(E) \geq g\left ( \tfrac 1 2 \arccos(t_{p,d}) \right )^2 \geq g(\arccos(t_{p,d}))^2 2^{-2d} = p^2 2^{-2d} \geq p^2 2^{-\tfrac{1}{2} \log(1/p) } \geq  p^{3} \left( 1 + c \left( \log \tfrac{1}{p} \right)^{3/2} \right)
\]
for a universal constant $c>0$, yielding the desired estimates. 

\medskip

From this point on, we may therefore assume that $d \geq \tfrac 1 4 \log \left( 1/p \right)$. We first prove the lemma under the assumption that $p < \tfrac 1 2$. The proof of the case $p \geq \tfrac 1 2$ is similar, as explained below. 

\medskip 

Let $Z$ be the first coordinate of a uniform point in $\mathbb{S}^{d-1}$, hence a random variable whose probability density function is $f_{d}(\cdot)$. Denote by $E_{i, j}$ the event $\{\langle X_i, X_j \rangle \geq t_{p, d}\}$ and by $E_{i, j} (x)$ the event $\{\langle X_i, X_j \rangle =x\}$. In what follows, we occasionally allow ourselves to condition on the zero probability event $E_{i,j}(x)$. This should be understood as conditioning on $\{\langle X_i, X_j \rangle \in [x-\epsilon,x+\epsilon] \}$ and then taking $\epsilon \to 0$, which is well-defined thanks to a simple continuity argument. 

\medskip 

 Note that $\P\left(E_{1,3}, E_{2, 3} \, \middle| \, E_{1, 2} (x) \right)$ is an increasing function of $x$. Therefore, 
\begin{align} \label{eq-trivial-1}
\P& \left(E_{1, 3}, E_{2, 3} \, \middle| \, E_{1, 2}) - \P(E_{1, 3}, E_{2, 3} \right) \nonumber \\
& =  \P \left(E_{1, 3}, E_{2, 3} \, \middle| \, E_{1, 2} \right) - \tfrac{1}{2}\P \left(E_{1, 3}, E_{2, 3} \, \middle| \, \langle X_1, X_2 \rangle \leq 0 \right) - \tfrac{1}{2}\P \left( E_{1, 3}, E_{2, 3} \, \middle| \, \langle X_1, X_2 \rangle \geq 0 \right) \nonumber \\
& \geq \tfrac{1}{2}  \left\{ \P \left( E_{1,3}, E_{2, 3} \, \middle| \, E_{1, 2} (t_{p, d}) \right) - \P \left( E_{1, 3}, E_{2, 3} \, \middle| \, E_{1, 2} (0) \right) \right\} \,.
\end{align}
Write $Z_1 = \langle X_3, X_1\rangle$ and $Z_2 = \left \langle X_3, \frac{\Proj_{X_1^\perp} X_2}{\left|\Proj_{X_1^\perp} X_2 \right|} \right \rangle$ where $\Proj_{X_1^\perp}$ denotes the orthogonal projection onto the orthogonal complement of $X_1$. We have that 
\begin{align} \label{eq-trivial-2}
&\P\left( E_{1,3}, E_{2, 3}\, \middle| \, E_{1, 2} (t_{p, d}) \right) - \P \left( E_{1, 3}, E_{2, 3} \, \middle| \, E_{1, 2} (0) \right) \nonumber \\
& \quad = \int_{z_1\geq t_{p, d}}  \left\{ \P  \left( E_{2, 3} \, \middle| \, Z_1 = z_1, E_{1, 2}(t_{p, d}) \right)  - \P   \left( E_{2, 3} \, \middle| \, Z_1 = z_1, E_{1, 2}(0) \right) \right\} f_{d}(z_1) d z_1 . 
\end{align}
Conditioning on $Z_1 = z_1$, it is easy to verify that $Z_2$ has the same distribution as $\sqrt{1-z_1^2} Z'$, where $Z'$ has density $f_{d-1}$.  Therefore 
\begin{align*} 
 \P  \left  ( E_{2, 3} \, \middle| \, Z_1 = z_1, E_{1, 2} (t_{p, d}) \right ) 
  &= \P \left (\sqrt{1-z_1^2}\sqrt{1- t_{p, d}^2} Z' + t_{p, d} z_1 \geq t_{p, d} \right ) \\
  & = \P \left (\sqrt{1- t_{p, d}^2} Z' \geq \sqrt{ \frac{1 - z_1}{1 + z_1}} t_{p, d} \right ).
\end{align*} 
Since $\sqrt{\left( 1- z_1 \right) / \left( 1 + z_1 \right)}$ is a decreasing function of $z_1$, the right hand side of the preceding inequality is increasing in $z_1$. 
Note also that $\P  \left(E_{2, 3} \, \middle| \, Z_1 = z_1, E_{1, 2}(0) \right)$ is a decreasing function of $z_1$. 
This is because under $E_{1,2} \left( 0 \right)$ we have $\left\langle X_2, X_3 \right\rangle = Z_2$, which, as mentioned above, under this conditioning has the same distribution as $\sqrt{1-z_1^2} Z'$, where $Z'$ has density $f_{d-1}$, and so the tail probability of $\left\langle X_2, X_3 \right\rangle$ is a decreasing function of $z_1$. 
Thus the integral in~\eqref{eq-trivial-2} is bounded from below by $p$ times the difference of the conditional probabilities evaluated at $z_1 = t_{p, d}$. 
Therefore, using also the previous display and that $\left\{ Z_1 = z_1 \right\} = E_{1,3} \left( z_1 \right)$, we have 
\begin{multline} \label{eqestd1}
\P \left(E_{1,3}, E_{2, 3}\, \middle| \, E_{1, 2} (t_{p, d}) \right) - \P\left(E_{1, 3}, E_{2, 3} \, \middle| \, E_{1, 2} (0) \right) \\
 \geq p \P \left ( \tfrac{t_{p,d}}{1 + t_{p,d}} \leq Z'\leq \tfrac{t_{p, d}}{\sqrt{1 - t_{p, d}^2}} \right ) \geq p \int_{\tfrac{t_{p,d}}{1+ t_{p,d}}}^{t_{p,d}} f_{d-1}(x) dx. 
\end{multline}
Since $f_{d-1}(z)/f_d(z) \geq 1/10$ for all $-1\leq z\leq 1$, using \eqref{eq:fdtpd} we learn that,
$$f_{d-1} (z) \geq c' d p t_{p, d} \mbox{ for all } 0\leq z\leq t_{p, d}\,,$$
where $c'>0$ is an absolute constant. Plugging the preceding inequality
into \eqref{eqestd1} gives that
\begin{align*}
\P \left(E_{1,3}, E_{2, 3}\, \middle| \, E_{1, 2} (t_{p, d}) \right) - \P\left(E_{1, 3}, E_{2, 3} \, \middle| \, E_{1, 2} (0) \right)  \geq c d p^2 t_{p, d}^3  \stackrel{\eqref{tpdest}}{\geq} c' (\tfrac{1}{2} - p)^3 p^2 \left (\tfrac{(\log (1/p))^{3/2}}{\sqrt{d}} \wedge d \right ) \,,
\end{align*}
where $c, c'>0$ are absolute constants and the last inequality follows from \eqref{tpdest}. By our assumption that $d \geq \tfrac{1}{4} \log \left( 1/p \right)$, we get that 
$$\P \left(E_{1,3}, E_{2, 3}\, \middle| \, E_{1, 2} (t_{p, d}) \right) - \P \left(E_{1, 3}, E_{2, 3} \, \middle| \, E_{1, 2} (0) \right)  \geq \frac{1}{16} c' (\tfrac{1}{2} - p)^3 p^2 \frac{\log (1/p)^{3/2}}{\sqrt{d}}\,,$$ 
Plugging the preceding inequality into \eqref{eq-trivial-1} we get that
\begin{align*}
\P \left( E \right) &= \P \left( E_{1,2}, E_{1,3}, E_{2,3} \right) = \P \left( E_{1,3}, E_{2,3} \, \middle| \, E_{1,2} \right) \P \left( E_{1,2} \right) \\ 
&= \left\{ \P \left( E_{1,3}, E_{2,3} \, \middle| \, E_{1,2} \right) - \P \left( E_{1,3}, E_{2,3}  \right) \right\} \P \left( E_{1,2} \right) + \P \left( E_{1,3}, E_{2,3} \right) \P \left( E_{1,2} \right) \\
&\geq \left\{ \frac{1}{32} c' \left( \tfrac{1}{2} - p \right)^3 p^2 \frac{\log (1/p)^{3/2}}{\sqrt{d}} \right\} \times p + p^2 \times p = p^3 \left( 1 + \frac{1}{32} c' \left( \tfrac{1}{2} - p \right)^3 \frac{\log (1/p)^{3/2}}{\sqrt{d}} \right),
\end{align*} 
which completes the proof of~\eqref{trianpsmall}, 
and it also completes the proof of~\eqref{triang} for the case $p<1/2$ (note that $\P(E_{1, 3}, E_{2, 3} ) = p^2$).

\medskip

It remains to prove the lemma for the case $p \geq 1/2$. Since the case $p=1/2$ is established in Lemma \ref{lem:tri_diff} in Section~\ref{sec:estimation}, we may assume that $\tfrac 1 2 < p < 1$. This case is treated by considering the event 
$\tilde E = E_{1, 2}^C \cap E_{1, 3}^C \cap E_{2, 3}^C $ and by observing that, since $\P(E_{1, 2}) = p$ and $\P(E_{1, 2}\cap E_{1, 3}) = p^2$, we have by the inclusion-exclusion principle that
\begin{equation}\label{eq-inc-exc}
\P(E) = 1 - 3p + 3p^2 - \P(\tilde E)\,.
\end{equation}
In light of this equation, we need to establish an upper bound on $\P(\tilde E)$. Our computation in what follows can be considered as a dual of the case for $p<1/2$. Note that $\P \left(E^C_{1,3}, E^C_{2, 3}\, \middle| \, E_{1, 2} (x) \right)$ is an increasing function of $x$. Therefore, 
\begin{align*}
\P \left(E_{1, 3}^C, E^C_{2, 3} \, \middle| \, E^C_{1, 2} \right) - \P(E^C_{1, 3}, E^C_{2, 3} ) \leq \tfrac{1}{2} \left(\P \left(E^C_{1,3}, E^C_{2, 3}\, \middle| \, E_{1, 2} (t_{p, d}) \right) - \P\left(E^C_{1, 3}, E^C_{2, 3} \, \middle| \, E_{1, 2} (0) \right) \right).
\end{align*}
We claim that the proof of the bound \eqref{triang} (and thus of the entire lemma) is concluded if we manage to show that
\begin{equation} \label{eq:ntslem1}
\P \left(E^C_{1,3}, E^C_{2, 3}\, \middle| \, E_{1, 2} (t_{p, d}) \right) - \P \left(E^C_{1, 3}, E^C_{2, 3} \, \middle| \, E_{1, 2} (0) \right)  \leq -c_p/\sqrt{d}\,,
\end{equation}
for a constant $c_p>0$ which depends only on $p$. Indeed, using the last two inequalities we would then have that 
\[
\P(\tilde E) = (1-p) \P \left(E_{1, 3}^C, E^C_{2, 3} \, \middle| \, E^C_{1, 2} \right) \leq - (1-p) c_p/\sqrt{d} + (1-p)^3,
\]
which, combined with \eqref{eq-inc-exc}, gives that $\P(E) \geq p^3 + (1-p) c_p/\sqrt{d}$. It thus remains to prove \eqref{eq:ntslem1}. 

\medskip

As above, we write $Z_1 = \langle X_3, X_1\rangle$ and $Z_2 = \left \langle X_3, \frac{\Proj_{X_1^\perp} X_2}{\left|\Proj_{X_1^\perp} X_2 \right|} \right \rangle$. We see that 
\begin{multline} \label{eq-trivial-2-1}
\P \left(E^C_{1,3}, E^C_{2, 3}\, \middle| \, E_{1, 2} (t_{p, d}) \right) - \P \left(E^C_{1, 3}, E^C_{2, 3} \, \middle| \, E_{1, 2} (0) \right) \\
= \int_{z_1\leq t_{p, d}}  \left( \P  \left(E^C_{2, 3} \, \middle| \, Z_1 = z_1, E_{1, 2}(t_{p, d}) \right)  - \P  \left(E^C_{2, 3} \, \middle| \, Z_1 = z_1, E_{1, 2}(0) \right) \right) f_{d}(z_1) d z_1 . 
\end{multline}
Conditioning on $Z_1 = z_1$, we see that $Z_2$ has the same distribution as $\sqrt{1-z_1^2} Z'$, where $Z'$ has density $f_{d-1}$.  Therefore 
\begin{align*} 
 \P  \left (E^C_{2, 3} \, \middle| \, Z_1 = z_1, E_{1, 2} (t_{p, d}) \right ) 
  &= \P \left (\sqrt{1-z_1^2}\sqrt{1- t_{p, d}^2} Z' + t_{p, d} z_1 \leq t_{p, d} \right ) \\
  & = \P \left (\sqrt{1- t_{p, d}^2} Z' \leq \sqrt{ \frac{1 - z_1}{1 + z_1}} t_{p, d} \right ).
\end{align*} 
Since $\sqrt{\left( 1- z_1 \right) / \left( 1 + z_1 \right)}$ is a decreasing function of $z_1$ (and $t_{p, d} <0$ in this case), the right hand side of the preceding inequality is increasing in $z_1$. It is clear that $ \P  \left(E^C_{2, 3} \, \middle| \, Z_1 = z_1, E_{1, 2}(0) \right)$ is a decreasing function of $z_1$. Thus, the right hand side of \eqref{eq-trivial-2-1} is bounded from above by $p$ times the integrand evaluated at $z_1 = t_{p, d}$. Therefore, 
\begin{align*}
\P \left(E^C_{1,3}, E^C_{2, 3}\, \middle| \, E_{1, 2} (t_{p, d}) \right) - \P \left(E^C_{1, 3}, E^C_{2, 3} \, \middle| \, E_{1, 2} (0) \right)  &\leq -p \P \left ( \tfrac{t_{p,d}}{1 + t_{p,d}} \leq Z'\leq \tfrac{t_{p, d}}{\sqrt{1 - t_{p, d}^2}} \right ) \\ 
&\leq -p t_{p,d} \left ( \tfrac{1}{\sqrt{1 - t_{p, d}^2}} -  \tfrac{1}{1+t_{p,d}} \right ) f_{d-1} \left ( \tfrac{t_{p,d}}{1+t_{p,d}} \right ).
\end{align*}
It is easily seen that the above is smaller than $-C_{p,d}$ where $C_{p,d}>0$ for any fixed $d$ and $\tfrac 1 2 < p < 1$. Thus, by choosing the constant $c_p$ to be small enough, we can make sure that the inequality in \eqref{eq:ntslem1} holds true for any fixed value of $d$. In other words, we may legitimately assume that $d > \tilde c_p$ where $\tilde c_p$ depends on $p$. Also, since $t_{1-p,d} = -t_{p,d}$, the second inequality in \eqref{tpdest} teaches us that
\begin{equation} \label{tpdsmall}
|t_{p,d}| \leq C' \sqrt{\frac{\log \tfrac{1}{1-p}}{d}},
\end{equation}
for a universal constant $C'>0$. In turn by taking $\tilde c_p$ to be large enough, we may assume that $|t_{p,d}| < 1/10$. Plugging this assumption to the last inequality, we deduce that
$$
\P \left(E^C_{1,3}, E^C_{2, 3}\, \middle| \, E_{1, 2} (t_{p, d}) \right) - \P \left(E^C_{1, 3}, E^C_{2, 3} \, \middle| \, E_{1, 2} (0) \right)  \leq -c p t_{p,d}^2 f_{d-1}(2 t_{p,d}).
$$
for a universal constant $c>0$. Using the bounds \eqref{gammaineq} and \eqref{tpdsmall}, we have
\[
f_{d-1} (2 t_{p,d})  \geq \tfrac{1}{2} \sqrt{d} (1-4t_{p,d}^2 )^{(d-3)/2} \geq \tfrac{1}{2} \sqrt{d} \left (1-4 (C')^2 \frac{\log \tfrac{1}{1-p}}{d} \right )^{(d-3)/2} \geq \sqrt{d} c_p'
\]
where $c_p'$ depends only on $p$. Combining the last two inequalities with~\eqref{tpdest} finally yields~\eqref{eq:ntslem1}. 
\end{proof}

\section{Proof of Theorem \ref{th:denseUB}} \label{secdenseUB}

Recall from the introduction that $T \left( G \right)$ and $\tau \left( G \right)$ denote the number of triangles and signed triangles, respectively, of a graph $G$ (see~\eqref{eq:triangle_def} and~\eqref{eq:signed_triangle_def}). 
We denote by $A$ the adjacency matrix of a graph $G$; we omit the dependence of $A$ on $G$, as the graph $G$ will always be obvious from the context. 
To abbreviate notation, for distinct vertices $i$, $j$, and $k$, we define $\overline{A}_{i,j} := A_{i,j} - \E \left[ A_{i,j} \right]$ and $\tau_G \left( i,j,k \right) := \overline{A}_{i,j} \overline{A}_{i,k} \overline{A}_{j,k}$. 
The number of signed triangles of a graph $G$ with vertex set $\left[ n \right]$ can then be written as follows:
\[
 \tau \left( G \right) = \sum_{\left\{ i, j, k \right\} \subset \binom{\left[n\right]}{3}} \tau_G \left( i, j, k \right).
\]

\medskip

In order to prove Theorem~\ref{th:denseUB}, we need to show that $\tau(G(n,p))$ and $\tau(G(n,p,d))$ behave very differently as long as $d$ is much smaller than $n^3$. 
This is done by establishing that the expectations of $\tau(G(n,p))$ and $\tau(G(n,p,d))$ differ by a quantity which exceeds both of the corresponding standard deviations. 
We break the proof into four parts accordingly: 
in Section~\ref{sec:tau_gnp} we analyze $\tau \left( G \left( n, p \right) \right)$, 
in Section~\ref{sec:tau_gnpd_exp} we give a lower bound on $\E \ \tau \left( G \left( n, p, d \right) \right)$, 
in Section~\ref{sec:tau_gnpd_var} we give an upper bound on the variance of $\tau \left( G \left( n, p, d \right) \right)$, 
and finally in Section~\ref{sec:pf_conc} we conclude the proof of Theorem~\ref{th:denseUB}. 

\subsection{Analysis of $\tau(G(n,p))$}\label{sec:tau_gnp}

The analysis of the statistic $\tau$ in the Erd\H{o}s-R\'enyi model is straightforward. 
First, note that for each $\left\{i,j\right\} \in  \binom{[n]}{2}$ one has $\EE\left[\overline{A}_{i,j}\right] = 0$. Therefore by the independence of edges it is clear that $\E \ \tau_{G(n,p)}(i,j,k) = 0$ for all $\left\{ i, j, k \right\} \subset \binom{\left[n\right]}{3}$, and so $\E \ \tau \left( G \left( n, p \right) \right) = 0$.

\medskip
 
In order to compute the variance, we again use the independence of edges to see that
\[
 \E \left[ \tau_{G\left(n,p \right)} \left( 1, 2, 3 \right) \tau_{G\left(n,p \right)} \left( 1, 2, 4 \right) \right] = \E \left[ \overline{A}_{1,2}^2 \overline{A}_{1,3} \overline{A}_{2,3} \overline{A}_{1,4} \overline{A}_{2,4} \right] = \E \left[ \overline{A}_{1,2}^2  \overline{A}_{1,3}  \overline{A}_{2,3}   \overline{A}_{1,4} \right] \E \left[  \overline{A}_{2,4} \right] = 0.
\]
Moreover, a simple calculation reveals that $\E \left[ \overline{A}_{i,j}^2 \right] = p \left( 1 - p \right)$ for every $\left\{i,j\right\} \in  \binom{[n]}{2}$. 
We conclude that for two triplets of vertices with indices  $i,j,k$ and $i',j',k'$ we have that
\[
  \E \left[ \tau_{G\left(n,p \right)} \left( i, j, k \right) \tau_{G\left(n,p \right)} \left( i', j', k' \right) \right] = 
\begin{cases}
0 & \text{if } \left\{i,j,k\right\} \neq \left\{i',j',k' \right\}, \\
p^3 \left( 1 - p \right)^3 & \text{if } \left\{i,j,k\right\} = \left\{i',j',k' \right\}.
\end{cases}
\]
Thus the variance of $\tau \left( G \left(n, p \right) \right)$ is
\[
\Var(\tau(G(n,p))) = \sum_{\left\{i,j,k \right\}} \sum_{\left\{i', j', k'\right\}} \E[\tau_{G(n,p)}(i,j,k) \tau_{G(n,p)}(i',j',k')] = \binom{n}{3} p^3 (1-p)^3.
\]
We have therefore established that
\begin{equation} \label{eq1-1}
\E \ \tau(G(n,p)) = 0 \;\; \text{and} \;\; \Var(\tau(G(n,p))) \leq \binom{n}{3}.
\end{equation}

\subsection{A lower bound for $\EE \ \tau(G(n,p,d))$}\label{sec:tau_gnpd_exp}

Our main goal here is to prove the following.

\begin{lemma} \label{lem:EtauGnpd}
For every $0 < p < 1$, there exists a constant $C_p > 0$ (depending only on $p$) such that for all $n$ and $d$ we have
\begin{equation} \label{eq1-2}
\E \ \tau(G(n, p, d)) \geq \binom{n}{3} \frac{C_p}{\sqrt{d}}.
\end{equation}
\end{lemma}
\begin{proof}
It suffices to estimate $\E \ \tau_{G\left( n,p,d \right)} \left( 1, 2, 3 \right)$, since $\E \ \tau(G(n, p, d)) = \binom{n}{3} \E \ \tau_{G\left( n,p,d \right)} \left( 1, 2, 3 \right)$. 
We have that
\begin{align*}
 \E \ \tau_{G\left( n,p,d \right)} \left( 1, 2, 3 \right) &= \E \overline{A}_{1,2} \overline{A}_{1,3} \overline{A}_{2,3} = \E \left( A_{1,2} - p \right) \left( A_{1,3} - p \right) \left( A_{2,3} - p \right) \\
&= \E A_{1,2} A_{1,3} A_{2,3} - p \left( \E A_{1,2} A_{1,3} + \E A_{1,2} A_{2,3} + \E A_{1,3} A_{2,3} \right) \\
&\qquad + p^2 \left( \E A_{1,2} + \E A_{1,3} + \E A_{2,3} \right) - p^3 \\
&= \E A_{1,2} A_{1,3} A_{2,3} - p^3,
\end{align*}
where the last equality follows from the simple facts that $\E A_{i,j} = p$ and $\E A_{i,j} A_{i,k} = p^2$ for all triples $\left\{ i, j, k \right\} \subset \binom{\left[n\right]}{3}$. 
The bound~\eqref{triang} from Lemma~\ref{lem:probatriangle-p} then gives that
\[
 \E \ \tau_{G\left( n,p,d \right)} \left( 1, 2, 3 \right) \geq \frac{C_p}{\sqrt{d}}
\]
for some $C_p > 0$, which concludes the proof. 
\end{proof}

\subsection{The variance of $\tau(G(n,p,d))$}\label{sec:tau_gnpd_var}

The estimation of the variance of $\tau(G(n,p,d))$ requires the following technical lemma.
\begin{lemma} \label{lem:var_tau_Gnpd}
For every $p \in \left[0,1 \right]$ we have that
\[
\E[ \tau_{G(n,p,d)}(1,2,3) \tau_{G(n,p,d)}(1,2,4)] \leq \pi^2/d.
\]
\end{lemma}

\begin{proof}
Define\footnote{The function $g$ in this proof should not be confused with the function $g$ appearing in the proof of Lemma~\ref{lem:probatriangle-p}.}
\[
g(x) := \E \left[ \tau_{G(n,p,d)}(1,2,3) \tau_{G(n,p,d)}(1,2,4) \, \middle| \, \langle X_1, X_2 \rangle = x \right].
\]
Note first that by rotational invariance, and by the independence of $X_3$ and $X_4$, we have that 
\begin{align*}
 g(x) &= \E\left[ \tau_{G(n,p,d)}(1,2,3) \tau_{G(n,p,d)}(1,2,4) \, \middle| \, X_1 = e_1, ~ X_2 = x e_1 + \sqrt{1-x^2} e_2 \right] \\
&= \left ( \mathbf{1}_{ \{x < t_{p,d} \} } p^2 + \mathbf{1}_{ \{x \geq t_{p,d} \} } \left(1 - p \right)^2 \right ) \left( \E \left[ \overline{A}_{1,3} \overline{A}_{2,3} \, \middle| \, X_1 = e_1, ~ X_2 = x e_1 + \sqrt{1-x^2} e_2 \right] \right)^2,
\end{align*}
where, as before, $\{e_1,\dots,e_d\}$ denotes the standard basis of $\R^d$. Define the spherical cap
\[
B_x := \left \{y \in \mathbb{S}^{d-1}: ~ \langle y, x e_1 + \sqrt{1-x^2} e_2 \rangle \geq t_{p,d} \right \}
\]
and let $S_x := B_1 \cap B_x$ denote the intersection of two such caps. 
Let $\sigma$ denote the normalized Haar measure on the sphere and note that $\sigma \left( B_1 \right) = \sigma \left( B_x \right) = p$. 
We have by definition that
\begin{multline*}
 \E \left[ \overline{A}_{1,3} \overline{A}_{2,3} \, \middle| \, X_1 = e_1, ~ X_2 = x e_1 + \sqrt{1-x^2} e_2 \right] \\
\begin{aligned}
&= \left( 1 - p \right)^2 \sigma \left( S_x \right) - p \left( 1 - p \right) \left( \sigma \left( B_1 \setminus S_x \right) +  \sigma \left( B_x \setminus S_x \right) \right) + p^2 \sigma \left( \mathbb{S}^{d-1} \setminus \left( B_1 \cup B_x \right) \right) \\
&= \left( 1 - p \right)^2 \sigma \left( S_x \right) - 2p \left( 1 - p \right) \left( p - \sigma \left( S_x \right) \right) + p^2 \left( 1 - 2p + \sigma \left( S_x \right) \right) = \sigma \left( S_x \right) - p^2,
\end{aligned}
\end{multline*}
and so
\[
 g \left( x \right) = \left ( \mathbf{1}_{ \{x < t_{p,d} \} } p^2 + \mathbf{1}_{ \{x \geq t_{p,d} \} } \left(1 - p \right)^2 \right ) \left( \sigma \left( S_x \right) - p^2 \right)^2 \leq \left( \sigma \left( S_x \right) - p^2 \right)^2 .
\]
Let $X$ be a random variable with the same law as that of $\langle X_1, X_2 \rangle$; in other words, $X$ has density $f_d$, see~\eqref{eq-density-Z} in Section~\ref{sec:triangle}. 
Thus by conditioning on $\left\langle X_1, X_2 \right\rangle$ and using the definition of $g\left( \cdot \right)$ we have that 
\begin{equation} \label{eqsigma}
\E \left[ \tau_{G(n,p,d)}(1,2,3) \tau_{G(n,p,d)}(1,2,4) \right] = \EE \left[g(X) \right] \leq \E \left[ \left( \sigma \left(S_X\right) - p^2 \right)^2 \right].
\end{equation}
Next, we claim that 
\begin{equation} \label{ESX}
\E \left[ \sigma(S_X) \right] = p^2.
\end{equation}
To see this, define, using a slight abuse of notation, $B_{\theta} = \left\{ y \in \mathbb{S}^{d-1} : \left\langle y, \theta \right \rangle \geq t_{p,d} \right\}$ for $\theta \in \mathbb{S}^{d-1}$, and let $Y$, $Z$, and $W$ be independent, uniformly distributed points on $\mathbb{S}^{d-1}$. By rotational invariance, $\sigma \left( S_X \right)$ has the same distribution as $\sigma \left( B_Y \cap B_Z \right)$, and so
\[
 \E \left[ \sigma(S_X) \right] = \E \left[ \sigma \left( B_Y \cap B_Z \right) \right] = \E \left[ \P \left(  W \in B_Y \cap B_Z \right) \right] = \P \left(  W \in B_Y \cap B_Z \right) = \left( \P \left( W \in B_Y \right) \right)^2 = p^2 
\]
as claimed. Thus by~\eqref{eqsigma} and~\eqref{ESX} we have that
\[
\E\left[ \tau_{G(n,p,d)}(1,2,3) \tau_{G(n,p,d)}(1,2,4) \right] \leq \Var \left[ h(X) \right],
\]
where we define $h(x) := \sigma \left( S_x \right)$. 

\medskip

Now since $1 = \E \left[ \left| Y \right|^2 \right] = \E \left[ Y_1^2 + \dots + Y_d^2 \right] = d \times \E \left[ Y_1^2 \right] = d \times \E \left[ X^2 \right]$ by linearity of expectation, 
where $Y = \left( Y_1, \dots, Y_d \right)$, we have that 
\[
 \EE \left[ X^2 \right] = \frac{1}{d}.
\]
We claim that the proof follows by showing that for all $-1 \leq x,y \leq 1$,
\begin{equation} \label{eq:Lipschitz}
|h(x) - h(y)| \leq 2 \left| \arcsin(x) - \arcsin(y) \right|.
\end{equation}
Indeed, observe that, since $\arcsin(x)$ is odd and convex on $[0,1]$, we have that
\begin{equation} \label{eq:convex}
\left| \arcsin(x) \right| \leq |x| \arcsin(1) = |x|\tfrac \pi 2, ~~ \forall -1 \leq x \leq 1,
\end{equation}
and now combining the last four displays we get that
\begin{align*}
\E\left[ \tau_{G(n,p,d)}(1,2,3) \tau_{G(n,p,d)}(1,2,4) \right] &\leq  \Var[h(X)] \leq \EE \left[ (h(X) - h(0))^2 \right] \\
&\stackrel{\eqref{eq:Lipschitz}}{\leq }  4 \EE \left[ (\arcsin(X))^2 \right] \stackrel{\eqref{eq:convex}}{\leq }  \pi^2 \EE \left[X^2 \right] = \tfrac{\pi^2}{d}, 
\end{align*}
which concludes the proof, using~\eqref{eq:Lipschitz}.

\medskip

It thus remains to prove inequality~\eqref{eq:Lipschitz}. To this end, consider the foliation of the sphere
\[
\mathbb{S}^{d-1} = \bigcup_{z \in \mathbf{B}^{d-2}} W_z,
\]
where $\mathbf{B}^{d-2}$ is the $(d-2)$-dimensional Euclidean unit ball and $W_z = \{ (x ,y , z); ~ x^2 + y^2 = 1 - |z|^2 \}$, with a corresponding decomposition $\theta \mapsto (x_\theta, y_\theta, z_\theta)$. Let $\mu$ denote the measure on $\mathbf{B}^{d-2}$ which is the push-forward of $\sigma$ under the map $\theta \mapsto z_\theta$, and for $z \in \mathbf{B}^{d-2}$ let $\sigma_{z}$ denote the uniform measure on $W_z$ so that 
\[
\int_{\mathbb{S}^{d-1}} f(\theta) d \sigma(\theta) = \int_{\mathbf{B}^{d-2}} \int_{W_z} f(x, y,  z) d\sigma_{z}(x,y) d \mu(z)
\]
for all measurable $f:\mathbb{S}^{d-1} \to \RR$. Clearly for all $0 < x,y < 1$ we have that
$$
\sigma(S_{x}) - \sigma(S_{y}) = \int_{\mathbf{B}^{d-2}} \left (\sigma_{z} (S_x \cap W_z) - \sigma_{z} (S_y \cap W_z) \right ) d \mu(z).
$$
By the triangle inequality, we have for all $z \in \mathbf{B}^{d-2}$ that 
$$
|\sigma_{z} (S_x \cap W_z) - \sigma_{z} (S_y \cap W_z) | \leq \sigma_z ((S_x \Delta S_y) \cap W_z) \leq \sigma_z ((B_x \Delta B_y) \cap W_z).
$$
Now fix $z \in \mathbf{B}^{d-2}$ and note that $W_z$ is a circle of radius $\sqrt{1 - \left| z \right|^2}$. 
If $\sqrt{1 - \left| z \right|^2} < t_{p,d}$, then $B_x \cap W_z = \emptyset$ for every $x \in \left[-1,1\right]$, and so $\sigma_z ((B_x \Delta B_y) \cap W_z) = 0$. 
If $\sqrt{1 - \left| z \right|^2} \geq t_{p,d}$, then $B_x \cap W_z$ and $B_y \cap W_z$ are both arcs of the circle $W_z$, 
and they can be transformed into each other by a rotation of angle $|\arcsin(x) - \arcsin(y)|$ in the appropriate direction. 
Consequently we have 
\[
\sigma_z ((B_x \Delta B_y) \cap W_z) \leq 2 |\arcsin(x) - \arcsin(y)| \sigma_z(W_z).
\]
Therefore, by the displays above and the triangle inequality, we have that
$$
|\sigma(S_{x}) - \sigma(S_{y})| \leq 2 |\arcsin(x) - \arcsin(y) | \sigma(\mathbb{S}^{d-1}) = 2 |\arcsin(x) - \arcsin(y) |,
$$
which is exactly \eqref{eq:Lipschitz}. 
\end{proof}

Using this lemma we can finally give an upper bound on the variance of $\tau(G(n,p,d))$. By expanding the variance as
\begin{multline} \label{sumvar}
\Var(\tau(G(n,p,d))) \\
= \sum_{\left\{i,j,k \right\}} \sum_{\left\{i', j', k' \right\}} 
\left\{ \E \left[ \tau_{G(n,p,d)}(i,j,k) \tau_{G(n,p,d)}(i',j',k') \right] - \E \left[ \tau_{G(n,p,d)}(i,j,k) \right] \E \left[ \tau_{G(n,p,d)}(i',j',k') \right] \right\}
\end{multline}
and by using the fact that $X_1,\dots,X_n$ are i.i.d., it is clearly enough to analyze the following terms:
\begin{align*}
W_1 &:= \E \left[ \tau_{G(n,p,d)}(1,2,3) \tau_{G(n,p,d)}(1,2,3) \right] - \E \left[ \tau_{G(n,p,d)}(1,2,3) \right] \E \left[ \tau_{G(n,p,d)}(1,2,3) \right], \\
W_2 &:= \E \left[ \tau_{G(n,p,d)}(1,2,3) \tau_{G(n,p,d)}(1,2,4) \right] - \E \left[ \tau_{G(n,p,d)}(1,2,3) \right] \E \left[ \tau_{G(n,p,d)}(1,2,4) \right], \\
W_3 &:= \E \left[ \tau_{G(n,p,d)}(1,2,3) \tau_{G(n,p,d)}(1,4,5) \right] - \E \left[ \tau_{G(n,p,d)}(1,2,3) \right] \E \left[ \tau_{G(n,p,d)}(1,4,5) \right], \\ 
W_4 &:= \E \left[ \tau_{G(n,p,d)}(1,2,3) \tau_{G(n,p,d)}(4,5,6) \right] - \E \left[ \tau_{G(n,p,d)}(1,2,3) \right] \E \left[ \tau_{G(n,p,d)}(4,5,6) \right].
\end{align*}
By symmetry, each term in the right hand side of~\eqref{sumvar} corresponds to one of these four ``types''. 
By definition we have $\left| W_1 \right| \leq 1$, and Lemma~\ref{lem:var_tau_Gnpd} states that $W_2 \leq \pi^2/d$. 
Since $\tau_{G(n,p,d)}(1,2,3)$ depends on $X_1$, $X_2$, and $X_3$, while 
$\tau_{G(n,p,d)}(4,5,6)$ depends on $X_4$, $X_5$, and $X_6$, 
by independence we have that $W_4 = 0$. 
To show that $W_3 = 0$, we can condition on $X_1$ and use rotational invariance: 
\begin{align*}
  \E \left[ \tau_{G(n,p,d)}(1,2,3) \tau_{G(n,p,d)}(1,4,5) \right] 
  &= \E \left[ \E \left[ \tau_{G(n,p,d)}(1,2,3) \tau_{G(n,p,d)}(1,4,5) \, \middle| \, X_1 \right] \right] \\
  &= \E \left[ \E \left[ \tau_{G(n,p,d)}(1,2,3) \, \middle| \, X_1 \right] \E \left[ \tau_{G(n,p,d)}(1,4,5) \, \middle| \, X_1 \right] \right] \\
  &= \E \left[ \E \left[ \tau_{G(n,p,d)}(1,2,3) \right] \E \left[ \tau_{G(n,p,d)}(1,4,5) \right] \right] \\
  &= \E \left[ \tau_{G(n,p,d)}(1,2,3) \right] \E \left[ \tau_{G(n,p,d)}(1,4,5) \right].
\end{align*}
The third equality above is due to the rotational invariance of $\left\{ X_i \right\}_{i=1}^n$, 
which implies that the conditional expectation $\E \left[ \tau_{G(n,p,d)}(1,2,3) \, \middle| \, X_1 \right]$ does not depend on $X_1$. 
Putting everything together, we finally conclude that 
\begin{equation}  \label{eq1-3}
\Var(\tau(G(n,p,d))) \leq \binom{n}{3} + \binom{n}{4} \binom{4}{2} \frac{\pi^2}{d} \leq n^3 + \frac{3 n^4}{d}.
\end{equation}

\subsection{Concluding the proof}\label{sec:pf_conc}

Combining displays~\eqref{eq1-1},~\eqref{eq1-2} and~\eqref{eq1-3} together, we get that
\[
\EE(\tau(G(n,p))) = 0, ~~ \EE[ \tau(G(n,p,d)) ] \geq \binom{n}{3} \frac{C_p}{\sqrt{d}}
\]
and
\[
\max \left\{ \Var\left[\tau(G(n,p)) \right], \Var\left[\tau(G(n,p,d)) \right] \right\} \leq n^3 + \frac{3 n^4}{d}.
\]
Using Chebyshev's inequality this implies that
\[
 \P\left(\tau(G(n,p,d)) \leq \tfrac{1}{2} \EE[ \tau(G(n,p,d)) ]  \right) \leq 200 \frac{d n^3 + 3 n^4  }{ C_p^2 n^6 },
\]
and also
$$\P\left(\tau(G(n,p)) \geq \tfrac{1}{2} \EE\left[ \tau(G(n,p,d)) \right] \right) \leq 200 \frac{d n^3 + 3 n^4  }{ C_p^2 n^6 } .$$
Putting the two above displays together we thus have that
$$\TV\left(\tau \left( G(n,p,d) \right), \tau \left( G(n,p) \right) \right) \geq 1 - \frac{Cd}{n^3} - \frac{C}{n^2} $$
for a constant $C$ that depends only on $p$. This concludes the proof of Theorem~\ref{th:denseUB}.

\section{Proof of Theorem \ref{sparseUB}}
Fix a constant $c > 0$ and let $p=c / n$.
To abbreviate notation, for distinct vertices $i$, $j$, and $k$ of a graph $G$, 
let $T_{G} \left( i, j , k \right) := A_{i,j} A_{i,k} A_{j,k}$; in other words, $T_{G} \left( i,j, k \right)$ is the indicator that the vertices $i$, $j$, and $k$ form a triangle in $G$. 
We can then write the number of triangles in a graph $G$ with vertex set $\left[n \right]$ as
\[
 T \left( G \right) = \sum_{\left\{i,j,k \right\} \subset \binom{\left[n\right]}{3}} T_G \left( i, j, k \right).
\]

\medskip

First, we have a simple estimate for the expectation of $T \left( G \left( n, c/n \right) \right)$:
\begin{equation} \label{eq3-1}
 \E \left[ T \left( G \left( n, c/n \right) \right) \right] = \binom{n}{3} \left( \frac{c}{n} \right)^3 \leq c^3.
\end{equation}

\medskip

We now turn to estimating the first two moments of $G(n,p,d)$. 
An application of the inequality~\eqref{trianpsmall} from Lemma~\ref{lem:probatriangle-p} gives that
\begin{equation} \label{eq3-3}
 \E \left[ T \left( G \left( n, c/n, d \right) \right) \right] \geq \kappa c^3 \frac{\left( \log \tfrac{n}{c} \right)^{3/2}}{\sqrt{d}}
\end{equation}
for a universal constant $\kappa > 0$. Note that the right hand side of the inequality~\eqref{eq3-3} goes to infinity when $d / \log^3 \left( n \right) \to 0$. 

\medskip

In order to establish an upper bound for the variance of $T \left( G \left( n, c/n, d \right) \right)$, first note that $T_{G\left( n, p, d \right)} \left( i, j, k \right)$ and $T_{G\left( n, p, d \right)} \left( i', j', k' \right)$ are independent whenever $\left\{i,j,k\right\}$ and $\left\{i',j',k'\right\}$ do not share an edge, i.e., whenever $\left|\left\{i,j,k\right\} \cap \left\{i',j',k'\right\} \right| \leq 1$.  
Furthermore, using the independence of $T_{G\left( n, p, d \right)} \left( 1, 2, 3 \right)$ and $A_{1,4}$, we have that
\[
 \E \left[ T_{G\left( n, p, d \right)} \left( 1, 2, 3 \right) T_{G\left( n, p, d \right)} \left( 1, 2, 4 \right) \right] \leq \E \left[ T_{G\left( n, p, d \right)} \left( 1, 2, 3 \right) A_{1,4} \right] = \E \left[ T_{G\left( n, p, d \right)} \left( 1, 2, 3 \right)  \right] \times p.
\]
Using these facts and expanding $T \left( G \left( n, p , d \right) \right)^2$ as a sum of indicators, we have that
\begin{multline*}
 \E \left[ T \left( G \left( n, p , d \right) \right)^2 \right] = \sum_{\left\{i,j,k\right\}} \sum_{\left\{i',j',k' \right\}} \E \left[ T_{G\left( n, p, d \right)} \left( i, j, k \right) T_{G\left( n, p, d \right)} \left( i', j', k' \right) \right] \\
\begin{aligned}
 &\leq \binom{n}{3}^2 \E \left[ T_{G\left( n, p, d \right)} \left( 1,2,3 \right) T_{G\left( n, p, d \right)} \left( 4,5,6 \right) \right] \\
&\qquad + \binom{n}{4} \binom{4}{2} \E \left[ T_{G\left( n, p, d \right)} \left( 1,2,3 \right) T_{G\left( n, p, d \right)} \left( 1,2,4 \right) \right] + \binom{n}{3} \E \left[ T_{G\left( n, p, d \right)} \left( 1,2,3 \right)^2 \right] \\
&\leq \left( \E \left[ T \left( G \left( n, p, d \right) \right) \right] \right)^2 + \binom{n}{3} \left( 2 n p + 1 \right) \E \left[ T_{G\left( n, p, d \right)} \left( 1,2,3 \right) \right],
\end{aligned}
\end{multline*}
and so
\begin{equation}\label{eq3-4}
 \Var \left[ T \left( G \left( n, c/n , d \right) \right) \right] \leq \left( 2c + 1 \right) \E \left[ T \left( G\left( n, c/n, d \right) \right) \right].
\end{equation}

\medskip

Now using this estimate, together with Chebyshev's inequality, gives that
\[
 \P \left( T \left( G \left( n, \tfrac{c}{n}, d \right) \right) \leq \tfrac{1}{2} \E \left[ T \left( G \left( n, \tfrac{c}{n}, d \right) \right) \right] \right) \leq \frac{4 \Var \left( T \left( G \left( n, \tfrac{c}{n}, d \right) \right) \right)}{\left( \E \left[ T \left( G \left( n, \tfrac{c}{n}, d \right) \right) \right] \right)^2} \leq \frac{8c+4}{\E \left[ T \left( G \left( n, \tfrac{c}{n}, d \right) \right) \right]}.
\]
On the other hand, Markov's inequality and the estimate~\eqref{eq3-1} together give that
\[
 \P \left( T \left( G \left( n, \tfrac{c}{n} \right) \right) \geq \tfrac{1}{2} \E \left[ T \left( G \left( n, \tfrac{c}{n}, d \right) \right) \right] \right) \leq \frac{2 c^3}{\E \left[ T \left( G \left( n, \tfrac{c}{n}, d \right) \right) \right]}.
\]
Finally, note that both of these upper bounds go to $0$ when $d / \log^3 \left( n \right) \to 0$, as can be seen from the estimate~\eqref{eq3-3}. 
This concludes the proof of Theorem~\ref{sparseUB}.

\section{Proof of the lower bound}\label{sec:denseLB}

Recall that if $Y_1$ is a standard normal random variable in $\R^d$, then $Y_1 / \left\| Y_1 \right\|$ is uniformly distributed on the sphere $\mathbb{S}^{d-1}$. 
Consequently we can view $G\left( n, p, d \right)$ as a function of an appropriate Wishart matrix, as follows. 
Let $Y$ be an $n \times d$ matrix where the entries are i.i.d.\ standard normal random variables, and let $W \equiv W (n,d) = YY^T$  be the corresponding $n \times n$ Wishart matrix. 
Note that $W_{ii} = \left\langle Y_i, Y_i \right\rangle = \left\| Y_i \right\|^2$ and so 
$\left\langle Y_i / \left\| Y_i \right\|, Y_j / \left\| Y_j \right\| \right\rangle = W_{ij} / \sqrt{W_{ii} W_{jj}}$. 
Thus the $n \times n$ matrix $A$ defined as
\[
 A_{i,j} =
\begin{cases}
 1 & \text{if } W_{ij} / \sqrt{W_{ii} W_{jj}} \geq t_{p,d} \text{ and } i \neq j\\
 0 & \text{otherwise}
\end{cases}
\]
has the same law as the adjacency matrix of $G\left(n,p,d\right)$. 
Denote the map that takes $W$ to $A$ by $H_{p,d}$, i.e., $A = H_{p,d} \left( W \right)$.

\medskip

In a similar way we can view $G \left( n, p \right)$ as a function of an $n \times n$ matrix drawn from the Gaussian Orthogonal Ensemble (GOE). 
Let $M\left( n \right)$ be a symmetric $n \times n$ random matrix where the diagonal entries are i.i.d.\ normal random variables with mean zero and variance 2, and the entries above the diagonal are i.i.d.\ standard normal random variables, with the entries on and above the diagonal all independent. 
Then the $n \times n$ matrix $B$ defined as 
\[
 B_{i,j} =
\begin{cases}
 1 & \text{if } M\left(n\right)_{i,j} \geq \overline{\Phi}^{-1} \left( p \right) \text{ and } i \neq j\\
 0 & \text{otherwise}
\end{cases}
\]
has the same law as the adjacency matrix of $G\left(n,p\right)$. 
Note that $B$ only depends on the off-diagonal elements of $M\left(n \right)$, so in the definition of $B$ we can replace $M\left( n \right)$ with $M \left( n, d \right) := \sqrt{d} M \left( n \right) + d I_n$, where $I_n$ is the $n \times n$ identity matrix, provided we also replace $\overline{\Phi}^{-1} \left( p \right)$ with $\sqrt{d} \times \overline{\Phi}^{-1} \left( p \right)$. 
Denote the map that takes $M\left(n,d\right)$ to $B$ by $K_{p,d}$, i.e., $B = K_{p,d} \left( M \left( n, d \right) \right)$. 
The maps $H_{p,d}$ and $K_{p,d}$ are different, but very similar, as we shall quantify later. 

\medskip

Using the triangle inequality, 
we can conclude from the previous two paragraphs that for any $p \in \left[0,1\right]$ we have that
\begin{multline*}
 \TV\left(G(n,p,d), G(n,p)\right) = \TV \left( H_{p,d} \left( W \left( n, d \right) \right), K_{p,d} \left( M \left( n,d \right) \right) \right) \\
\begin{aligned}
 &\leq \TV \left( H_{p,d} \left( W \left( n, d \right) \right), H_{p,d} \left( M \left( n,d \right) \right) \right) + \TV \left( H_{p,d} \left( M \left( n, d \right) \right), K_{p,d} \left( M \left( n,d \right) \right) \right) \\
 &\leq \TV \left(  W \left( n, d \right) ,  M \left( n,d \right) \right)  + \TV \left( H_{p,d} \left( M \left( n, d \right) \right), K_{p,d} \left( M \left( n,d \right) \right) \right). 
\end{aligned}
\end{multline*}
The second term in the line above is a smaller order error term which we deal with in Appendix~\ref{sec:error}, 
and thus 
Theorem~\ref{th:mainresult}\eqref{th:main_denseLB} 
follows from the following result, which is a more precise restatement of Theorem~\ref{thm:Wishart_GOE_intro} from the Introduction.
\begin{theorem}\label{thm:Wishart_GOE}
 Define the random matrix ensembles $W\left( n, d \right)$ and $M \left( n, d \right)$ as above. If $d / n^3 \to \infty$, then
\begin{equation}\label{eq:Wishart_GOE}
  \TV \left( W \left( n, d \right), M \left( n, d \right) \right) \to 0.
\end{equation}
\end{theorem}

\medskip

After proving this result, we learned of the work of Jiang and Li~\cite{jiang2013approximation}, who also study the question of when can Wishart matrices be approximated by a GOE random matrix. 
They prove a result analogous to Theorem~\ref{thm:Wishart_GOE} for the joint eigenvalue distributions, from which Theorem~\ref{thm:Wishart_GOE} follows by orthogonal invariance. 
Nonetheless, we present here our proof, for the following two reasons:
(1)~to keep the paper self-contained, and 
(2)~because our proof is shorter and simpler than the one presented in~\cite{jiang2013approximation}. 

\medskip

Theorem~\ref{thm:Wishart_GOE} thus states that as $d/n^3 \to \infty$, all statistics of the random matrix ensembles $W\left(n,d \right)$ and $M\left(n, d \right)$ have asymptotically the same distribution. 
We note that in the random matrix literature there has been lots of work showing that particular statistics (e.g., the empirical spectral distribution, or the largest eigenvalue) of these two ensembles have asymptotically the same distribution, even when $d / n \to c \in \left[ 0, \infty \right]$, see, e.g.,~\cite{bai2009spectral,johnstone2001distribution,elkaroui2003largest,karoui2007tracy}.

\medskip

\begin{proof}
We show~\eqref{eq:Wishart_GOE} by comparing the densities of the two random matrix ensembles. Let $\mathcal{P} \subset \R^{n^2}$ denote the cone of positive semidefinite matrices. 
It is well known (see, e.g.,~\cite{wishart1928generalised}) that when $d \geq n$, $W(n,d)$ has the following density with respect to the Lebesgue measure on $\mathcal{P}$:
\[
 f_{n,d} \left( A \right) := \frac{\left( \det \left( A \right) \right)^{\frac{1}{2} \left( d - n - 1 \right)} \exp \left( - \frac{1}{2} \Tr \left( A \right) \right)}{2^{\frac{1}{2}dn} \pi^{\frac{1}{4} n \left( n-1\right)} \prod_{i=1}^n \Gamma \left( \frac{1}{2} \left( d+1-i \right) \right)},
\]
where $\Tr\left( A \right)$ denotes the trace of the matrix $A$. 
It is also known that the density of a GOE random matrix with respect to the Lebesgue measure on $\R^{n^2}$ is
$A \mapsto \left( 2 \pi \right)^{-\frac{1}{4} n \left( n + 1 \right)} 2^{-\frac{n}{2}} \exp \left( - \frac{1}{4} \Tr \left( A^2 \right) \right)$ 
and so the density of $M \left( n, d \right)$ with respect to the Lebesgue measure on $\R^{n^2}$ is
\[
 g_{n,d} \left( A \right) := \frac{\exp \left( - \frac{1}{4d} \Tr \left( \left( A - d I_n \right)^2 \right) \right)}{\left( 2\pi d \right)^{\frac{1}{4} n \left( n + 1 \right)} 2^{\frac{n}{2}}}.
\]
Denote the measure given by this density by $\mu_{n,d}$, let $\lambda$ denote the Lebesgue measure on $\R^{n^2}$ and write $A \succeq 0$ if $A$ is positive semidefinite.
We can then write
\begin{align*}
 \TV \left( W(n,d), M \left( n, d \right) \right) &= \frac{1}{2} \int_{\R^{n^2}} \left| f_{n,d} \left( A \right) \mathbf{1}_{\left\{ A \succeq 0 \right\}} - g_{n,d} \left( A \right) \right| d \lambda \left( A \right) \\
 &= \frac{1}{2} \int_{\R^{n^2}} \left| \frac{f_{n,d} \left( A \right) \mathbf{1}_{\left\{ A \succeq 0 \right\}}}{g_{n,d} \left( A \right)} - 1 \right| d \mu_{n,d} \left( A \right).
\end{align*}
In order to show that $\TV \left( W(n,d), M \left( n, d \right) \right) = o \left( 1 \right)$, it suffices to show that 
\begin{equation} \label{eq:likelihood_ratio}
 \frac{f_{n,d} \left( A \right) \mathbf{1}_{\left\{ A \succeq 0 \right\}}}{g_{n,d} \left( A \right)} = 1 + o \left( 1 \right) 
\end{equation} 
as $d/n^3 \to \infty$ with probability $1- o \left( 1 \right)$ according to the measure $\mu_{n,d}$. 
This is because
\[
 \int_{\R^{n^2}} \left( \frac{f_{n,d} \left( A \right)\mathbf{1}_{\left\{ A \succeq 0 \right\}}}{g_{n,d} \left( A \right)} - 1 \right) d \mu_{n,d} \left( A \right) = 0,
\]
and so 
\begin{equation}\label{eq:pos_neg}
 \int_{\R^{n^2}} \left( \frac{f_{n,d} \left( A \right)\mathbf{1}_{\left\{ A \succeq 0 \right\}}}{g_{n,d} \left( A \right)} - 1 \right)_{+} d \mu_{n,d} \left( A \right)
 = \int_{\R^{n^2}} \left( \frac{f_{n,d} \left( A \right)\mathbf{1}_{\left\{ A \succeq 0 \right\}}}{g_{n,d} \left( A \right)} - 1 \right)_{-} d \mu_{n,d} \left( A \right).
\end{equation}
Since 
$\left| \left( \frac{f_{n,d} \left( A \right)\mathbf{1}_{\left\{ A \succeq 0 \right\}}}{g_{n,d} \left( A \right)} - 1 \right)_{-} \right| \leq 1$, 
if~\eqref{eq:likelihood_ratio} holds, then the right hand side of~\eqref{eq:pos_neg} is $o(1)$, and consequently so is the left hand side of~\eqref{eq:pos_neg}, which then shows that $\TV \left( W(n,d), M \left( n, d \right) \right) = o \left( 1 \right)$. 

\medskip

It is known (see, e.g.,~\cite{anderson2010introduction}) that, with probability $1-o\left(1\right)$, all the eigenvalues of $M(n)$ are in the interval $\left[-3 \sqrt{n}, 3 \sqrt{n} \right]$ and so all the eigenvalues of  $M \left( n, d \right)$ are in the interval $\left[d - 3 \sqrt{dn}, d + 3 \sqrt{dn} \right]$. 
Since $d/n^3 \to \infty$, we thus have $\P \left( M \left( n, d \right) \succeq 0 \right) = 1 - o \left( 1 \right)$, 
and so we may restrict our attention to $\mathcal{P}$. 
Define $\alpha_{n,d} \left( A \right) := \log \left( f_{n,d} \left( A \right) / g_{n,d} \left( A \right) \right)$. 
It remains then to show that $\alpha_{n,d} \left( A \right) = o \left( 1 \right)$ as $d / n^3 \to \infty$ with probability $1 - o \left( 1 \right)$ according to the measure $\mu_{n,d}$. 

\medskip

Denote the eigenvalues of an $n \times n$ matrix $A$ by $\lambda_1 \left( A \right) \leq \dots \leq \lambda_n \left( A \right)$; when the matrix is obvious from the context, we omit the dependence on $A$. Recall that $\det \left(A \right) = \prod_{i=1}^n \lambda_i$ and $\Tr \left( A \right) = \sum_{i=1}^n \lambda_i$. We then have
\begin{align*}
 \alpha_{n,d} \left( A \right) &= \frac{1}{2} \sum_{i=1}^n \left\{ \left( d - n - 1 \right) \log \lambda_i - \lambda_i + \frac{1}{2d} \left( \lambda_i - d \right)^2 \right\} \\
 &\quad + \left\{ \frac{n\left( n + 3 \right)}{4} - \frac{dn}{2} \right\} \log 2 + \frac{n}{2} \log \pi + \frac{n\left( n + 1 \right)}{4} \log d  - \sum_{i=1}^n \log \Gamma \left( \frac{1}{2} \left( d + 1 - i \right) \right).
\end{align*}
By Stirling's formula we know that $\log \Gamma \left( z \right) = \left( z - \frac{1}{2} \right) \log z - z + \frac{1}{2} \log \left( 2 \pi \right) + O \left( \frac{1}{z} \right)$ as $z \to \infty$, so 
\begin{align*}
 \alpha_{n,d} \left( A \right) &= \frac{1}{2} \sum_{i=1}^n \left\{ \left( d - n - 1 \right) \log \lambda_i - \lambda_i + \frac{1}{2d} \left( \lambda_i - d \right)^2 \right\} \\
 &\quad + \frac{n\left( n + 1 \right)}{4} \log d  - \frac{1}{2} \sum_{i=1}^n \left( d - i \right) \log \left( d + 1 - i \right) + \frac{1}{2} \sum_{i=1}^n \left( d + 1 - i \right) + O \left( \frac{n}{d} \right).
\end{align*}
Now writing $\log \left( d + 1 - i \right) = \log d + \log \left( 1 - \frac{i-1}{d} \right) = \log d  - \frac{i-1}{d} + O \left( \frac{i^2}{d^2} \right)$ we get that
\begin{align*}
 \alpha_{n,d} \left( A \right) &= \frac{1}{2} \sum_{i=1}^n \left\{ \left( d - n - 1 \right) \log \lambda_i - \lambda_i + \frac{1}{2d} \left( \lambda_i - d \right)^2 \right\} \\
 &\quad + \frac{1}{2} \left\{ - n d \log d + nd + n \left( n + 1 \right)  \log d \right\} + O \left( \frac{n^3}{d} \right).
\end{align*}
Defining $h \left( x \right) := \frac{1}{2} \left\{ \left( d - n - 1 \right) \log \left( x / d \right) - \left( x - d \right) + \frac{1}{2d} \left( x - d \right)^2 \right\}$, we have that
\begin{equation}\label{eq:alpha_sum_f}
 \alpha_{n,d} \left( A \right)  = \sum_{i=1}^n h\left( \lambda_i \right) + O \left( \frac{n^3}{d} \right).
\end{equation}
Note that the derivatives of $h$ at $d$ are $h \left( d \right) = 0$, $h'\left( d \right) = - \frac{n+1}{2d}$, $h'' \left( d \right) = \frac{n+1}{2d^2}$, and $h^{\left( 3 \right)} \left( d \right) = \frac{d - n - 1}{d^3}$, and also $h^{\left( 4 \right)} \left( x \right) = - \frac{3 \left( d - n - 1 \right)}{x^4}$. Approximating $h$ with its third order Taylor polynomial around $d$ we get
\[
 h(x) = - \frac{n+1}{2d} \left( x - d \right) + \frac{n+1}{4d^2} \left( x - d \right)^2 + \frac{d - n - 1}{6d^3} \left( x - d \right)^3 - \frac{d-n-1}{8 \xi^4} \left( x - d \right)^4,
\]
where $\xi$ is some real number between $x$ and $d$.
By~\eqref{eq:alpha_sum_f} we need to show that $\sum_{i=1}^n h\left( \lambda_i \right) = o \left( 1 \right)$ with probability tending to $1$. 
Recall that with probability $1 - o \left( 1 \right)$ all eigenvalues of  $M \left( n, d \right)$ are in the interval $\left[d - 3 \sqrt{dn}, d + 3 \sqrt{dn} \right]$. If $\lambda_i \in \left[d - 3 \sqrt{dn}, d + 3 \sqrt{dn} \right]$  and $\xi_i$ is between $\lambda_i$ and $d$ for all $i \in \left[n \right]$, then it is immediate that 
\[
 \sum_{i=1}^n \left| \frac{n+1}{4d^2} \left( \lambda_i - d \right)^2 - \frac{d-n-1}{8 \xi_i^4} \left( \lambda_i - d \right)^4 \right| \leq \frac{c n^3}{d}
\]
for some constant $c$, and so what remains to show is that
\begin{equation}\label{eq:eig_symmetry}
 \sum_{i=1}^n \left\{ - \frac{n+1}{2d} \left( \lambda_i - d \right) + \frac{d - n - 1}{6d^3} \left( \lambda_i - d \right)^3 \right\} = o \left( 1 \right)
\end{equation}
as $d / n^3 \to \infty$ with probability $1 - o \left( 1 \right)$. This follows from known results about the moments of the empirical spectral distribution of Wigner matrices. 
In particular,~\cite[Theorem~2.1.31]{anderson2010introduction} shows that 
\[
 \frac{1}{\left( n d \right)^{1/2}} \sum_{i=1}^n \left( \lambda_i - d \right) \qquad \text{ and } \qquad \frac{1}{\left( n d \right)^{3/2}} \sum_{i=1}^n \left( \lambda_i - d \right)^3
\]
both converge weakly to a normal distribution with appropriate variance, from which~\eqref{eq:eig_symmetry} immediately follows.
\end{proof}

\section{Dimension estimation} \label{sec:estimation}

In this section we prove Theorem~\ref{thm:estimation}. 
The idea for the proof is very similar to that of Theorem~\ref{th:denseUB}, 
and also uses the statistic $\tau(G)$ counting the ``number'' of signed triangles, analyzed in Section~\ref{secdenseUB}. 
However, dimension estimation is a slightly more subtle matter, 
since here it is necessary to have a bound on the difference of the expected number of triangles between consecutive dimensions, 
rather than just a lower bound on the expected number of triangles in the random geometric graph $G\left(n,p,d\right)$, 
as in Lemma~\ref{lem:probatriangle-p}. 
The next lemma gives a bound on this difference; but note that this lemma only deals with the case $p=1/2$. 
We believe that this result should hold true for any fixed $0<p<1$, but the proof seems to be much more involved.

\begin{lemma} \label{lem:tri_diff}
Let $\left\{e_i : 1 \leq i \leq d \right\}$ denote the standard basis in $\mathbb{R}^d$, and define, 
for each dimension $d$,  
\[
h(d) := \P \left(\langle X_1, e_1 \rangle \geq 0, \langle X_2, e_1 \rangle \geq 0, \langle X_1, X_2 \rangle \geq 0 \right),
\]
where $X_1$ and $X_2$  are independent uniformly distributed random vectors on $\mathbb{S}^{d-1}$. Then
\begin{equation} \label{eqdifftriang}
h(d-1) - h(d) \geq c d^{-3/2}
\end{equation}
for some universal constant $c>0$.
\end{lemma}
\begin{proof}
Define the event $E := \{ \langle X_1, e_1 \rangle \geq 0, \langle X_2, e_1 \rangle \geq 0, \langle X_1, X_2 \rangle \geq 0 \}$. 
We first claim that  
\begin{equation} \label{wedge}
\P \left(E \, \middle| \, \langle X_1, e_1 \rangle =  t \right) = \frac{\pi - \arccos t}{2 \pi} \mathbf{1}_{\left\{ t \in \left[0,1\right] \right\}}
\end{equation}
for all $t \in \left[-1,1 \right]$. 
To see this, first note that by the rotational invariance of the vectors $X_1$ and $X_2$, and since the event $E$ is invariant under the action of the orthogonal group on the three vectors simultaneously, 
when conditioning on $\left\langle X_1, e_1 \right\rangle = t$ 
we can assume that $X_1 =  t e_1 + \sqrt{1-t^2} e_2$. 
Thus for all $t \in \left[0,1\right]$ we have that 
\begin{equation} \label{eq:wedge}
\P \left( E \, \middle| \, \langle X_1, e_1 \rangle = t \right) = \P \left (\langle X_2, e_1 \rangle \geq 0, t \langle X_2, e_1 \rangle + \sqrt{1-t^2} \langle X_2, e_2 \rangle \geq 0  \right ).
\end{equation}
Since the projection of $X_2$ on the span of $e_1$ and $e_2$ is distributed according to a law which is invariant under rotations of this span, 
we have that the probability in~\eqref{eq:wedge} is proportional to the angle of the wedge defined by $\{x \in \mathrm{span}(e_1,e_2) : \langle x, e_1 \rangle \geq 0, t \langle x, e_1 \rangle + \sqrt{1-t^2} \langle x, e_2 \rangle \geq 0 \}$. The angle of this wedge is exactly $\pi - \arccos t$, and thus formula \eqref{wedge} follows. 

\medskip 

Let $Z_d$ be a random variable with the same law as that of $\left\langle X_1, e_1 \right\rangle$; in other words, $Z_d$ has density $f_d$, see~\eqref{eq-density-Z} in Section~\ref{sec:triangle}. Defining $Y_d = \arccos \left( Z_d \right)$ we therefore have that
\[
 \P \left( E \right) = \E \left[ \left( \frac{1}{2} - \frac{Y_d}{2\pi} \right) \mathbf{1}_{\left\{ Y_d \in \left[ 0, \pi/2 \right] \right\}} \right] = \frac{1}{4} - \frac{1}{2\pi} \E \left[ Y_d \mathbf{1}_{\left\{ Y_d \in \left[ 0, \pi/2 \right] \right\}} \right]
\]
Using the change of variables $y = \arccos z$ and the formula for $f_d$ in~\eqref{eq-density-Z}, we have that the density of $Y_d$ is given by
\[
 g_d \left( x \right) = f_d \left( \cos \left( x \right) \right) \sin \left( x \right) = \frac{\Gamma \left( d/2 \right)}{\Gamma \left( (d-1) / 2 \right) \sqrt{\pi}} \left( \sin \left( x \right) \right)^{d-2}, ~~ x \in \left[ 0, \pi \right].
\]
To abbreviate notation, define $A_d := \frac{\Gamma \left( d/2 \right)}{\Gamma \left( (d-1) / 2 \right) \sqrt{\pi}}$. 
By integrating we thus get that
\begin{equation} \label{eq:hd}
 h \left( d \right) = \P \left( E \right) = \frac{1}{4} - \frac{A_d}{2\pi} \int_0^{\pi/2} x \left( \sin \left( x \right) \right)^{d-2} dx.
\end{equation}
Let $\mu_d$ be the probability measure on $\left[0, \pi / 2 \right]$ defined by 
\[
 \frac{d \mu_d \left( x \right)}{dx} = 2 A_d \left( \sin \left( x \right) \right)^{d-2}.
\]
We then have that $\frac{d \mu_d}{d \mu_{d-1}} \left( x \right) = \frac{A_d}{A_{d-1}} \sin \left( x \right)$, and so 
\begin{equation} \label{eq:mu_identity}
 \int_{0}^{\pi/2} \sin \left( x \right) d \mu_{d-1} \left( x \right) = \frac{A_{d-1}}{A_d}.
\end{equation}
Elementary estimates concerning the $\Gamma$ function give that 
\begin{equation} \label{eq:gamma}
 \frac{A_{d}}{A_{d-1}} = \frac{\Gamma \left( \frac{d}{2} \right) \Gamma \left( \frac{d-2}{2} \right)}{\left( \Gamma \left( \frac{d-1}{2} \right) \right)^2} \in \left( 1 + \frac{1}{12d}, 1 + \frac{2}{d} \right), \qquad  \text{ and } \qquad \frac{A_{d-1}}{A_d} \in \left( 1 - \frac{2}{d}, 1 - \frac{1}{12d} \right).
\end{equation}
In the following we write $\nu \equiv \mu_{d-1}$ to abbreviate notation. We then have that
\begin{align}
 h(d-1) - h(d) &= \frac{1}{2\pi} \left( A_d \int_0^{\pi/2} x \left( \sin \left( x \right) \right)^{d-2} dx - A_{d-1} \int_0^{\pi/2} x \left( \sin \left( x \right) \right)^{d-3} dx  \right) \notag \\
&= \frac{1}{4\pi} \left( \int x d \mu_d \left( x \right) - \int x d \mu_{d-1} \left( x \right) \right) \notag \\
&= \frac{1}{4\pi} \int x \left( \frac{A_d}{A_{d-1}} \sin \left( x \right) - 1 \right) d \nu \left( x \right) \notag \\
&= \frac{1}{4\pi} \frac{A_d}{A_{d-1}} \int x \left( \sin \left( x \right) - \frac{A_{d-1}}{A_d} \right) d \nu \left( x \right) \notag \\
&= \frac{1}{4\pi} \frac{A_d}{A_{d-1}} \int x \left( \sin \left( x \right) - \int \sin \left( y \right) d\nu \left( y \right) \right) d\nu \left( x \right). \label{eq:hd-1-hd}
\end{align}
Since the function $x \mapsto \sin \left( x \right)$ is continuous and monotone on the interval $\left[0,\pi/2\right]$, there exists a unique $x_0 \in \left(0, \pi/2 \right)$ such that 
\begin{equation}\label{eq:x0}
 \sin \left( x_0 \right) = \int_0^{\pi/2} \sin \left( y \right) d \nu \left( y \right).
\end{equation}
Since the function $x \mapsto \sin \left( x \right)$ is concave on the interval $\left[0,\pi/2\right]$, 
we have that $\sin \left( x \right) - \sin \left( x_0 \right) \leq \cos \left( x_0 \right) \left( x - x_0 \right)$ for all $x \in \left[0,\pi/2\right]$. Thus continuing from~\eqref{eq:hd-1-hd} and also using that $A_{d} / A_{d-1} \geq 1$ from~\eqref{eq:gamma}, we have that
\begin{align}
 h\left( d- 1 \right) - h \left( d \right) &= \frac{1}{4\pi} \frac{A_d}{A_{d-1}} \int_0^{\pi/2} \left( x - x_0 \right) \left( \sin \left( x \right) - \sin \left( x_0 \right) \right) d\nu \left( x \right) \notag \\
 &\geq \frac{1}{4\pi} \cos \left( x_0 \right) \int_0^{x_0} \left( x - x_0 \right)^2 d \nu \left( x \right). \label{eq:estdiff1}
\end{align}
Using~\eqref{eq:mu_identity},~\eqref{eq:gamma}, and~\eqref{eq:x0}, we have that 
\[
 x_0 = \arcsin \left( \frac{A_{d-1}}{A_d} \right) \in \left( \arcsin \left( 1 - \frac{2}{d} \right), \arcsin \left( 1 - \frac{1}{12d} \right) \right),
\]
and so
\[
 \cos \left( x_0 \right) \geq \cos \left( \arcsin \left( 1 - \tfrac{1}{12d} \right) \right) = \cos \left( \arccos \left( \sqrt{1 - \left( 1 - \tfrac{1}{12d} \right)^2} \right) \right) = \sqrt{\frac{1}{6d} - \frac{1}{144d^2}} \geq \frac{1}{4 \sqrt{d}}.
\]
Furthermore, since $x \mapsto \arcsin \left( x \right)$ is convex on $\left[0,1\right]$ with derivative $1/ \sqrt{1-x^2}$, we have that for all $x \in \left[ 0, \arcsin \left( 1 - \tfrac{3}{d} \right) \right]$, 
\[
 \left( x - x_0 \right)^2 \geq \left( \arcsin \left( 1 - \tfrac{2}{d} \right) - \arcsin \left( 1 - \tfrac{3}{d} \right) \right)^2 \geq \left( \frac{1}{d} \times \frac{1}{\sqrt{1 - \left( 1 - \tfrac{3}{d} \right)^2}} \right)^2 = \frac{1}{6d - 9} \geq \frac{1}{6d}.
\]
Plugging the estimates from the previous three displays back into~\eqref{eq:estdiff1}, we get that
\begin{equation} \label{eq:estdiff2}
 h \left( d - 1 \right) - h \left( d \right) \geq \frac{1}{100 \pi d^{3/2}} \nu \left( \left[ 0, \arcsin \left( 1 - \tfrac{3}{d} \right) \right] \right).
\end{equation}
Let $Z'$ be a random variable with density $f_{d-1}$. Then by the definition of $\nu$ we have that 
\begin{align*}
 \nu \left( \left[ 0, \arcsin \left( 1 - \tfrac{3}{d} \right) \right] \right) &= 2 \P \left( \arccos \left( Z' \right) \in \left[ 0, \arcsin \left( 1 - \tfrac{3}{d} \right) \right] \right) = 2 \P \left( Z' \geq \cos \left( \arcsin \left( 1 - \tfrac{3}{d} \right) \right) \right) \\
&= 2 \P \left( Z' \geq \cos \left( \arccos \left( \sqrt{ 1 - \left( 1 - \tfrac{3}{d} \right)^2 } \right) \right) \right) \geq 2 \P \left( Z' \geq \tfrac{\sqrt{6}}{\sqrt{d-1}} \right), 
\end{align*}
which, by~\eqref{eq-Psi-lower}, is bounded from below by a universal constant. Together with~\eqref{eq:estdiff2}, this concludes the proof. 
\end{proof}

\begin{proof}\textbf{of Theorem \ref{thm:estimation}} 
Recall from the proof of Lemma~\ref{lem:EtauGnpd} that 
\[
 \E \left[ \tau \left( G\left( n, \tfrac{1}{2} , d \right) \right) \right] = \binom{n}{3} \left( h \left( d \right) - \left( 1/2 \right)^3 \right). 
\]
Thus by Lemma~\ref{lem:tri_diff} we have for all $d_1 < d_2$ that  
\begin{align*}
 \E \left[ \tau \left( G \left( n, \tfrac{1}{2}, d_1 \right) \right) \right] - \E \left[ \tau \left( G \left( n, \tfrac{1}{2}, d_2 \right) \right) \right] &= \binom{n}{3} \left( h \left( d_1 \right) - h \left( d_2 \right) \right) \geq \binom{n}{3} \left( h \left( d_1 \right) - h \left( d_1 + 1 \right) \right) \\
&\geq \binom{n}{3} \frac{c}{\left( d_1 + 1 \right)^{3/2}} \geq \frac{c_1 n^3}{d_1^{3/2}}
\end{align*}
for a universal constant $c_1 > 0$. By~\eqref{eq1-3} we have a bound on the variance of these statistics:
\[
 \max\left\{ \Var \left( \tau \left( G \left( n, \tfrac{1}{2}, d_1 \right) \right)  \right), \Var \left( \tau \left( G \left( n, \tfrac{1}{2}, d_2 \right) \right)  \right) \right\} \leq n^3 + \frac{3n^4}{d_1}.
\]
Using Chebyshev's inequality thus gives that
\[
 \P \left( \tau \left( G \left( n, \tfrac{1}{2}, d_1 \right) \right) \leq  \frac{1}{2} \left( \E \left[ \tau \left( G \left( n, \tfrac{1}{2}, d_1 \right) \right) \right] + \E \left[ \tau \left( G \left( n, \tfrac{1}{2}, d_2 \right) \right) \right] \right) \right) \leq \frac{12}{c_1^2} \frac{n^3 d_1^3 + n^4 d_1^2}{n^6} \leq \frac{24}{c_1^2} \frac{d_1^2}{n^2},
\]
where the second inequality holds if $d_1 \leq n$; if $d_1 > n$ then the upper bound on the probability is vacuously true by choosing $c_1 \leq 1$. Similarly we also have that
\[
 \P \left( \tau \left( G \left( n, \tfrac{1}{2}, d_2 \right) \right) \geq  \frac{1}{2} \left( \E \left[ \tau \left( G \left( n, \tfrac{1}{2}, d_1 \right) \right) \right] + \E \left[ \tau \left( G \left( n, \tfrac{1}{2}, d_2 \right) \right) \right] \right) \right) \leq \frac{12}{c_1^2} \frac{n^3 d_1^3 + n^4 d_1^2}{n^6} \leq \frac{24}{c_1^2} \frac{d_1^2}{n^2}.
\]
Putting the two previous displays together concludes the proof, and shows that we can take $C = 48\max\left\{1/ c_1^2,1\right\}$ for the constant in the statement of the theorem.
\end{proof}


\section*{Acknowledgements}

We thank Noureddine El Karoui, G\'abor Lugosi, and Johan Ugander for helpful discussions and useful references, and Benedek Valk\'o for pointing us to ref.~\cite{jiang2013approximation}. 
We also thank Shirshendu Ganguly and three anonymous referees for useful comments that helped improve the manuscript. 
This work was done while J.D.\ and R.E.\ were visiting researchers and M.Z.R.\ was an intern at the Theory Group of Microsoft Research. They thank the Theory Group for their hospitality. 
J.D. acknowledges support from NSF grant DMS 1313596, and M.Z.R.\  acknowledges support from NSF grant DMS 1106999.


\bibliographystyle{plain}
\bibliography{bib}


\appendix

\section{The error term in the proof of the lower bound}\label{sec:error}

Recalling notation from Section~\ref{sec:denseLB}, we prove here the following lemma, which is necessary in order to conclude the proof of Theorem~\ref{th:mainresult}\eqref{th:main_denseLB}. 
\begin{lemma}\label{lem:errorLB}
 If $d / n^3 \to \infty$, then $\sup_{p \in \left[0,1\right]} \TV \left( H_{p,d} \left( M \left( n, d \right) \right), K_{p,d} \left( M \left( n,d \right) \right) \right) \to 0$. 
\end{lemma}
Define $\Psi_d \left( x \right) = \int_{x / \sqrt{d}}^1 f_d \left( y \right) dy$, and note that by the definition of $t_{p,d}$, $p = \int_{t_{p,d}}^1 f_d \left( x \right) dx = \Psi_d \left( t_{p,d} \sqrt{d} \right)$. 
During the proof of Lemma~\ref{lem:errorLB} we use the following result of Sodin~\cite{Sodin05}.
\begin{lemma}\label{lem:sodin}
 There exist constants $C, C_1, C_2 > 0$, and a sequence $\epsilon_d \searrow 0$ that satisfies $\eps_d = O \left( 1 / d \right)$, such that the following inequalities hold for all $0 < t < C \sqrt{d}$:
\[
\left( 1 - \eps_d \right) \overline{\Phi} \left( t \right) \exp \left( - C_1 t^4 / d \right) \leq \Psi_d \left( t \right) \leq \left( 1 + \eps_d \right) \overline{\Phi} \left( t \right) \exp \left( -  C_2 t^4 / d \right).
\]
\end{lemma}
We note that in~\cite[Lemma~1]{Sodin05} the fact that $\eps_d$ satisfies $\eps_d = O \left( 1 / d \right)$ is not specified, but this can be read from the proof, where this error comes from the error in  Stirling's formula. 

\medskip

\begin{proof}\textbf{of Lemma~\ref{lem:errorLB}} 
To abbreviate notation, we write $X := H_{p,d} \left( M \left( n, d \right) \right)$, $Y := K_{p,d} \left( M \left( n, d \right) \right)$, and also $M \equiv M \left( n \right)$. 
Define $\wh{D}_{i,j} := \sqrt{\left( 1 + M_{i,i} / \sqrt{d} \right) \left( 1 + M_{j,j} / \sqrt{d} \right)}$ and recall that for $1 \leq i < j \leq n$ we have that
\begin{align}
 X_{i,j} &= \mathbf{1}_{\left\{ \wh{D}_{i,j}^{-1} M_{i,j}  \,  \geq \,  \sqrt{d} t_{p,d}  \right\}}, \label{eq:Xij}\\
 Y_{i,j} &= \mathbf{1}_{\left\{ M_{i,j} \,  \geq \,  \overline{\Phi}^{-1} \left( p \right) \right\}}, \notag
\end{align}
and define also
\begin{equation}\label{eq:Zij}
 Z_{i,j} = \mathbf{1}_{\left\{ M_{i,j} \geq \sqrt{d} t_{p,d} \right\}}.
\end{equation}
By the triangle inequality we have that
\[
 \TV \left( X, Y \right) \leq \TV \left( X, Z \right) + \TV \left( Y, Z \right),
\]
and we deal with the second term first. 
By a union bound, we have that
\[
 \TV \left( Y, Z \right) \leq \P \left( \exists \, 1 \leq i < j \leq n \, : \, Y_{i,j} \neq Z_{i,j} \right) \leq n^2 \P \left( Y_{1,2} \neq Z_{1,2} \right), 
\]
and thus it remains to show that $\P \left( Y_{1,2} \neq Z_{1,2} \right) = o \left( n^{-2} \right)$; since $n^3 / d \to 0$, it is enough to show that $\P \left( Y_{1,2} \neq Z_{1,2} \right) = O \left( d^{-2/3} \right)$. Since $M_{1,2}$ is a standard normal random variable, we have that 
\begin{equation}\label{eq:YZ}
 \P \left( Y_{1,2} \neq Z_{1,2} \right) = \left| \overline{\Phi} \left( \sqrt{d} t_{p,d} \right) - p \right| = \left| \overline{\Phi} \left( \sqrt{d} t_{p,d} \right) - \Psi_d \left( \sqrt{d} t_{p,d} \right) \right|.
\end{equation}
For $p = 1/2$ this expression is zero, and noting that $t_{p,d} = - t_{1-p,d}$ by symmetry, it is enough to bound the expression in~\eqref{eq:YZ} for $p \in [0,1/2)$. 
For $p \in \left( 0, 1 / 2 \right)$ fixed, note that Lemma~\ref{cor:tpd} states that there exists a constant $C_p < \infty$ such that $0 \leq \sqrt{d} t_{p,d} \leq C_p$. 
Together with Lemma~\ref{lem:sodin} applied at $t = \sqrt{d} t_{p,d}$, this implies that the expression in~\eqref{eq:YZ} is $O \left( 1 / d \right)$. 

\medskip 

More generally, 
by Lemma~\ref{cor:tpd} there exists a universal constant $C > 0$, such that if $p \in \left( n^{-\alpha}, 1/2 \right)$, then $0 \leq \sqrt{d} t_{p,d} \leq C \sqrt{\alpha \log \left( n \right) }$, and if $p \in \left[0, n^{-\alpha} \right]$, then $\sqrt{d} t_{p,d} \geq C^{-1} \sqrt{\alpha \log \left( n \right)}$. 
When $p \in \left( n^{-\alpha}, 1/2 \right)$, Lemma~\ref{lem:sodin} applied at $t = \sqrt{d} t_{p,d}$ then implies that the expression in~\eqref{eq:YZ} is $O \left( \alpha^2 \log^2 \left( n \right) / d \right)$, which is $O \left( d^{-2/3} \right)$ for constant $\alpha$. When $p \in \left[ 0, n^{-\alpha} \right]$ we then have that
\[
 \P \left( Y_{1,2} \neq Z_{1,2} \right) \leq p + \overline{\Phi} \left( \sqrt{d} t_{p,d} \right) \leq n^{-\alpha} + \overline{\Phi} \left( C^{-1} \sqrt{\alpha \log \left(n \right)} \right) \leq n^{-\alpha} + \exp \left( - \frac{\alpha \log\left(n \right)}{2C^2} \right).
\]
By choosing $\alpha := \max \left\{ 3, 6 C^2 \right\}$, this expression becomes $O\left( n^{-3} \right)$. This concludes the proof that 
\[
 \sup_{p \in \left[0,1\right]} \P \left( Y_{1,2} \neq Z_{1,2} \right) = O \left( d^{-2/3} + n^{-3} \right). 
\]

\medskip 

We now deal with the term $\TV \left( X, Z \right)$. 
For a matrix $A$, let $D' \left( A \right)$ denote the diagonal matrix obtained from $A$ by setting the non-diagonal terms to zero, and let $D \left( A \right) := A - D' \left( A \right)$. 
With $I \equiv I_n$ denoting the $n \times n$ identity matrix, define the function 
\[
f \left( A \right) := \left( I_n + \frac{1}{\sqrt{d}} D' \left( A \right) \right)^{-1/2} A \left( I_n + \frac{1}{\sqrt{d}} D' \left( A \right) \right)^{-1/2}.
\]
Note that $Z$ can be obtained from $D \left( M \right)$ by thresholding the entries appropriately (see~\eqref{eq:Zij}), and $X$ can be obtained from $D \left( f \left( M \right) \right)$ in the same way (see~\eqref{eq:Xij}). 
Consequently we have that
\begin{equation} \label{eq:TVproj}
\sup_{p \in \left[0,1\right]} \TV \left( X, Z \right) \leq \TV \left(  D \left( f \left( M \right) \right), D\left( M \right) \right),
\end{equation}
and in the remainder of the proof we show that this latter quantity goes to zero as $n^3 / d \to 0$. 

\medskip 
Let $\Omega = \RR^{(n^2 - n) / 2}$,  which we canonically identify with the space of symmetric $n \times n$ matrices with zero on the diagonal, 
and let $\Omega' = \RR^n$, which we think of as the space of the diagonal entries on an $n \times n$ matrix. By slight abuse of notation, we allow ourselves to think of $D$ and $D'$ as functions from the space of symmetric $n \times n$ matrices to $\Omega$ and $\Omega'$ respectively. We can thus naturally view the function $D \circ f$ as a mapping from $\Omega \oplus \Omega'$ to $\Omega$. 

\medskip 

In order  to bound the right hand side of \eqref{eq:TVproj}, we need an estimate for the density of $D \left(f \left(M \right) \right)$, which we denote by $w(x)$. Note that  $w(x)$ is the density of the push-forward under the function $D \circ f$ of the measure whose density is $\gamma(x,y)$, where
$$
\gamma(x,y) := \frac{1}{2^{n/2} (2 \pi)^{(n^2+n)/2}} \exp \left ( - \tfrac 1 2 \sum_{1 \leq i<j \leq n} x_{i,j}^2 - \tfrac 1 4 \sum_{i=1}^{n} y_i^2 \right )
$$
is the Gaussian GOE density on $\Omega \oplus \Omega'$. 
%
%
Eventually, we would like to compare $w(x)$ to the density of $D(M)$, which we denote by $q(x)$. Since $q(x)$ is the push-forward of $\gamma$ under $D$, we have that 
\begin{equation*}
q(x) = \frac{1}{(2 \pi)^{(n^2-n)/2}} \exp \left ( - \tfrac 1 2 \sum_{1 \leq i<j \leq n} x_{i,j}^2 \right )
\end{equation*}
which is just the standard Gaussian density on $\Omega$.

\medskip

Our next goal is thus to provide a formula for $w(x)$, which we do in a rather straightforward manner, by integration over fibers of the map $(D \circ f)^{-1}$. To this end, fix $x \in \Omega$. Observe that, by definition, $D'$ is a bijection from $(D \circ f)^{-1}(x)$ to $\Omega'$. In other words, for $y \in \Omega'$, defining
$$
g(x,y) := \left(\frac{1}{\sqrt{d}} y + I_n \right)^{1/2} x \left( \frac{1}{\sqrt{d}} y + I_n  \right)^{1/2},
$$
it can easily be checked by the definition of $f$ and $g$ that
$$
(D \circ f)^{-1} (x) = \{ (D(g(x,y)), y); ~ y \in \Omega' \}.
$$
Define a mapping $T:\Omega \oplus \Omega' \to \Omega \oplus \Omega'$ by
$$
T(x,y) := (D(g(x,y)), y).
$$
Let $\phi(\cdot)$ be a test function on $\Omega$. Then, substituting $(x,y) = T(x',y')$ and noting that $D(f(T(x',y'))) = x'$, we have that
\begin{multline*}
  \int_{\Omega} \phi(x) w(x) dx \\
\begin{aligned}
&= \int_{\Omega \times \Omega'} \phi(D(f(x,y))) \gamma(x,y) dxdy = \int_{\Omega \times \Omega'} \phi(D(f(T(x',y')))) \gamma(T(x',y')) J(x',y') dx' dy' \\
&= \int_{\Omega \times \Omega'} \phi(x') \gamma(T(x',y')) J(x',y') dx' dy' = \int_{\Omega} \phi(x') \left (\int_{\Omega'} \gamma(T(x',y')) J(x',y') dy' \right ) dx' 
\end{aligned}
\end{multline*}
where $J(\cdot,\cdot)$ is the Jacobian of $T$. Since the above holds true for any measurable test function, recalling the definition of $T$ we conclude that
\begin{equation} \label{eq:wx}
w(x) = \int_{\Omega'} \gamma(D(g(x,y)), y) J(x,y) dy.
\end{equation}
Let us now calculate the Jacobian $J(\cdot,\cdot)$. Recall that $T:\Omega \oplus \Omega' \to \Omega \oplus \Omega'$. Writing $\nabla T$ as a block matrix corresponding to this decomposition, we clearly have
\[
\nabla T = \left ( 
\begin{matrix}
\frac{\partial g}{\partial x} & \frac{\partial g}{\partial y} \\
0 & I_n
\end{matrix}
\right ).
\]
Therefore, we have $J(x,y) = \det\left ( \frac{\partial g}{\partial x} \right )$. Now by the definition of $g$ we have that 
$$
\frac{\partial g_{i,j}}{\partial x_{i',j'}} = \mathbf{1}_{\{i=i',j=j'\}} \left(\frac{1}{\sqrt{d}} y_i + 1 \right)^{1/2} \left(\frac{1}{\sqrt{d}} y_j + 1 \right)^{1/2},
$$
so we finally get that
\begin{equation*}
J(x,y) = \prod_{1 \leq i < j \leq n} \left|\frac{1}{\sqrt{d}} y_i + 1 \right|^{1/2} \left|\frac{1}{\sqrt{d}} y_j + 1 \right|^{1/2}  = \prod_{i=1}^n \left | 1 + \frac{y_i}{\sqrt{d}} \right |^{(n-1)/2}.
\end{equation*}
Together with the formula \eqref{eq:wx}, this yields
\begin{align*}
 w(x) &= \frac{1}{2^{n/2} (2 \pi)^{(n^2+n)/2}} \int_{\Omega'} J(x, y) \exp \left (- \frac 1 2 \sum_{1 \leq i<j \leq n} g(x,y)_{i,j}^2 - \frac 1 4 \sum_j y_j^2 \right ) dy\\
&= \frac{1}{2^{n/2} (2 \pi)^{(n^2+n)/2}} \int_{\Omega'} J(x, y) \exp \left (- \frac 1 2 \sum_{1 \leq i<j \leq n} x_{i,j}^2 \left (1 + \frac{y_i + y_j}{\sqrt{d}} + \frac{y_i y_j}{d} \right ) - \frac 1 4 \sum_j y_j^2 \right )dy \\
&= \frac{1}{(2 \pi)^{(n^2-n)/2}} \EE \left [ \prod_{i=1}^n \left | 1 + \frac{\Gamma_i}{\sqrt{d/2}} \right |^{(n-1)/2}  \exp \left ( - \frac 1 2 \sum_{1 \leq i<j \leq n} x_{i,j}^2 \left ( 1 + \frac{\Gamma_i + \Gamma_j }{ \sqrt{d/2}}  + \frac{2\Gamma_i \Gamma_j}{ d} \right ) \right ) \right ],
\end{align*}
where $\Gamma=(\Gamma_1,...,\Gamma_n)$ is a standard Gaussian random vector. We can now estimate
\begin{align*}
\frac{w(x)}{q(x)} 
&= \EE \left [ \prod_{i=1}^n \left | 1 + \frac{\Gamma_i}{\sqrt{d/2}} \right |^{(n-1)/2}  \exp \left ( - \frac 1 2 \sum_{1 \leq i<j \leq n} x_{i,j}^2 \left ( \frac{\Gamma_i + \Gamma_j}{ \sqrt{d/2}}  + \frac{2\Gamma_i \Gamma_j}{d} \right ) \right ) \right ] \\
&\geq \EE \left [ \prod_{i=1}^n \left | 1 + \frac{\Gamma_i}{\sqrt{d/2}} \right |^{(n-1)/2}  \exp \left ( - \frac 1 2 \sum_{1 \leq i<j \leq n} x_{i,j}^2 \left ( \frac{\Gamma_i + \Gamma_j}{ \sqrt{d/2}}  + \frac{2 \left( \Gamma_i^2 + \Gamma_j^2 \right)}{d} \right ) \right ) \right ] \\
&= \prod_{i=1}^n \EE \left [  \left | 1 + \frac{\Gamma_i}{\sqrt{d/2}} \right |^{(n-1)/2} \exp \left ( -  \sum_{j \neq i} x_{i,j}^2 \left ( \frac{\Gamma_i}{\sqrt{2d}} + \frac{\Gamma_i^2}{d} \right ) \right ) \right ] \\
&\geq \prod_i\EE \left [  \left( 1 + \frac{n}{2\sqrt{2d}} \Gamma_i \right) \left( 1 - \frac{\sum_{j \neq i}  x_{i,j}^2 }{\sqrt{2d}} \Gamma_i  - \frac{\sum_{j \neq i} x_{i,j}^2}{d} \Gamma_i^2 \right) \mathbf{1}_{\left\{ \left| \Gamma_i \right| \leq \sqrt{d} /2 \right\}} \right ] \\
&=  \prod_i\EE \left [ \left(  1 + \frac{\tfrac n 2 -\sum_{j \neq i}  x_{i,j}^2}{\sqrt{2d}} \Gamma_i - \frac{(n/4 + 1) \sum_{j \neq i} x_{i,j}^2 }{d} \Gamma_i^2 - \frac{n \sum_{j \neq i}  x_{i,j}^2}{2\sqrt{2}d^{3/2}} \Gamma_i^3  \right) \mathbf{1}_{\left\{ \left| \Gamma_i \right| \leq \sqrt{d} /2 \right\}} \right ] \\
&=  \prod_i\EE \left [ \left(  1  - \frac{(n/4 + 1) \sum_{j \neq i}  x_{i,j}^2 }{d} \Gamma_i^2   \right) \mathbf{1}_{\left\{ \left| \Gamma_i \right| \leq \sqrt{d} /2 \right\}} \right ].
\end{align*}
A well-known estimate concerning the Gaussian distribution gives $\P(|\Gamma_1| > \sqrt{d}/2 ) < 2 e^{-d/8}$. Consequently we have for all $1 \leq i \leq n$ that
\[
\EE \left [ \left(  1  - \frac{(n/4 + 1) \sum_{j \neq i} x_{i,j}^2 }{d} \Gamma_i^2   \right) \mathbf{1}_{\left\{ \left| \Gamma_i \right| \leq \sqrt{d} /2 \right\}} \right ] \geq 1 - 2 e^{-d/8} - \frac{(n/4 + 1) \sum_{j \neq i} x_{i,j}^2 }{d},
\]
and the two last displays give that 
\begin{align*}
\frac{w(x)}{q(x)} ~& \geq \left ( 1 - 2ne^{-d/8} - \frac{n/4 + 1}{d} \sum_{1 \leq i < j \leq n}   x_{i,j}^2 \right ) \mathbf{1}_{ \left  \{ \frac{(n/4 + 1) \sum_{j \neq i} x_{i,j}^2 }{d} < \tfrac 1 2, \forall 1 \leq i \leq n \right  \}  } \\
& \geq \left ( 1 - 2ne^{-d/8} - \frac{n/4 + 1}{d} \sum_{1 \leq i < j \leq n}   x_{i,j}^2 \right ) \mathbf{1}_{ \left  \{ x_{i,j}^2 < \tfrac{d}{n^2}, \forall 1 \leq i,j \leq n \right  \}  }.
\end{align*}
Finally, let $U=\{U_{i,j}\}_{1 \leq i < j \leq n}$ be a standard Gaussian vector on $\Omega$. The above display gives that 
\begin{align*}
  \TV \left( D \left( f \left( M \right) \right), D \left( M \right) \right) &= \E \left| \frac{w(D(M))}{q(D(M))} - 1 \right| = 2 \E \left( 1 - \frac{w(U)}{q(U)} \right)_{+} \\
  &\leq 2ne^{-d/8} + \frac{n/4 + 1}{d} \E \left[  \sum_{i,j} U_{i,j}^2 \right] + \P \left ( \max_{1\leq i < j \leq n} U_{i,j}^2 > \tfrac{d}{n^2} \right) \\
  &\leq 2ne^{-d/8} + \frac{2n^3}{d} + 2n^2 e^{- \tfrac{d}{2n^2}},
\end{align*}
which completes the proof.
\end{proof}

\end{document}